# CONSTRUCTIBILITÉ ET MODÉRATION UNIFORMES EN COHOMOLOGIE ÉTALE

Fabrice Orgogozo[‡]

Résumé. Soient $S$ un schéma nœthérien et $f : X → S$ un morphisme propre. D'après SGA 4 XIV, pour tout faisceau constructible $\mathscr{F}$ de $\mathbb{Z}/n\mathbb{Z}$-modules sur $X$, les faisceaux de $\mathbb{Z}/n\mathbb{Z}$-modules $\mathrm{R}^i f_\star \mathscr{F}$, obtenus par image directe (pour la topologie étale), sont également constructibles : il existe une stratification $\mathfrak{S}$ de $S$ telle que ces faisceaux soient localement constants constructibles sur les strates. À la suite de travaux de N. Katz et G. Laumon, ou L. Illusie, dans le cas particulier où $S$ est génériquement de caractéristique nulle ou bien les faisceaux $\mathscr{F}$ sont constants (de torsion inversible sur $S$), on étudie ici la dépendance de $\mathfrak{S}$ en $\mathscr{F}$. On montre qu'une condition naturelle de constructibilité et modération « uniforme » satisfaite par les faisceaux constants, introduite par O. Gabber, est stable par les foncteurs $\mathrm{R}^i f_\star$. Si $f$ n'est pas supposé propre, ce résultat subsiste sous réserve de modération à l'infini, relativement à $S$. On démontre aussi l'existence de bornes uniformes sur les nombres de Betti, qui s'appliquent notamment pour les fibres des faisceaux $\mathrm{R}^i f_\star \mathbb{F}_\ell$, où $\ell$ parcourt les nombres premiers inversibles sur $S$.

Abstract. Let $S$ be a Noetherian scheme and $f : X → S$ a proper morphism. By SGA 4 XIV, for any constructible sheaf $\mathscr{F}$ of $\mathbb{Z}/n\mathbb{Z}$-modules on $X$, the sheaves of $\mathbb{Z}/n\mathbb{Z}$-modules $\mathrm{R}^i f_\star \mathscr{F}$ obtained by direct image (for the étale topology) are themselves constructible, that is, there is a stratification $\mathfrak{S}$ of $S$ on whose strata these sheaves are locally constant constructible. After previous work of N. Katz and G. Laumon, or L. Illusie, on the special case in which $S$ is generically of characteristic zero or the sheaves $\mathscr{F}$ are constant (with invertible torsion on $S$), here we study the dependency of $\mathfrak{S}$ on $\mathscr{F}$. We show that a natural "uniform" tameness and constructibility condition satisfied by constant sheaves, which was introduced by O. Gabber, is stable under the functors $\mathrm{R}^i f_\star$. If $f$ is not proper, this result still holds assuming tameness at infinity, relative to $S$. We also prove the existence of uniform bounds on Betti numbers, in particular for the stalks of the sheaves $\mathrm{R}^i f_\star \mathbb{F}_\ell$, where $\ell$ ranges through all prime numbers invertible on $S$.

## Table des matières

[‡]CMLS, École polytechnique, CNRS, Université Paris-Saclay, 91128 Palaiseau Cedex, France
Fabrice.Orgogozo+math@normalesup.org

Introduction

**0.1.** Soient $S$ un schéma et $f : X \to S$ un morphisme propre de présentation finie. D'après un théorème de finitude classique, dû à Michael Artin et Alexander Grothendieck ([SGA 4 XIV, 1.1]), pour tout faisceau constructible $\mathscr{F}$ de $\mathbb{Z}/n\mathbb{Z}$-modules sur $X$, les faisceaux de $\mathbb{Z}/n\mathbb{Z}$-modules $\mathrm{R}^i f_\star \mathscr{F}$, obtenus par image directe, sont également *constructibles* : (lorsque $S$ est cohérent, par exemple nœthérien) il existe une stratification $\mathfrak{S}$ de $S$ telle que les faisceaux $\mathrm{R}^i f_\star \mathscr{F}$ soient localement constants et de présentation finie sur les strates. (Ils sont également nuls pour $i \gg 0$ indépendamment de $\mathscr{F}$.)

À $f$ fixé, il est naturel de s'interroger sur la dépendance de $\mathfrak{S}$ en $\mathscr{F}$, notamment sous la forme suivante, qui n'est qu'un cas particulier des questions que l'on souhaite aborder ici :

> Existe-t-il une stratification $\mathfrak{S}$ qui convienne pour tous les faisceaux $\mathrm{R}^i f_\star \mathbb{Z}/n\mathbb{Z}$, où $n$ parcourt les entiers inversibles sur $S$ ?

(Notons que si l'entier $n$ parcourt les puissances d'un même nombre premier $\ell$, comme cela est naturel du point de vue de la cohomologie $\ell$-adique, la réponse à la question précédente est trivialement « oui ».)

**0.2.** Si $S$ est de type fini sur le corps $\mathbb{C}$ des nombres complexes, l'existence d'une telle stratification est bien connue en vertu du résultat analogue pour le morphisme d'espaces topologiques $X(\mathbb{C}) \to S(\mathbb{C})$ — cf. par exemple [Borel 1984, V.10.16] ou [Beilinson, Bernstein et Deligne 1982, p. 153, premier paragraphe] — et du théorème de comparaison entre cohomologie de Betti et cohomologie étale de M. Artin et A. Grothendieck ([SGA 4 XVI, 4.1]).

**0.3.** Bien entendu, tant pour les applications que pour les démonstrations, il est utile de ne pas se restreindre à des faisceaux constants mais de considérer des faisceaux — et même des *complexes* — constructibles le long d'une stratification $\mathfrak{X}$ de $X$ (au sens évident rappelé en 1.1.3). Cependant, en cohomologie étale, des phénomènes de ramification sauvage font que l'on ne peut s'attendre à ce que, donnés un morphisme propre $f : X \to S$ et $\mathfrak{X}$ une stratification de $X$, il existe une stratification $\mathfrak{S}$ de $S$ telle que les images directes supérieures par $f$ de faisceaux sur $X$ constructibles le long de $\mathfrak{X}$ soient constructibles le long de $\mathfrak{S}$. En effet, si par exemple $\ell$ est un nombre premier fixé, $k$ un corps algébriquement clos de caractéristique $p > 0$ divisant $\ell - 1$ et $\psi : \mathbb{Z}/p\mathbb{Z} \to \mathbb{F}_\ell^\times$ un caractère non-trivial, les $\mathbb{F}_\ell$-faisceaux d'Artin-Schreier $\mathscr{L}_\lambda := \mathscr{L}_\psi((x - \lambda)y)$ sur $\mathbf{A}^1_{x/\mathbb{F}_p} \times \mathbf{A}^1_{y/\mathbb{F}_p}$ sont lisses pour chaque $\lambda \in k$ mais les faisceaux image directe exceptionnelle $\mathrm{R}^2 \mathrm{pr}_{1!} \mathscr{L}_\lambda$ ne sont pas constructibles le long d'une même stratification de $\mathbf{A}^1_{x/\mathbb{F}_p}$, indépendante de $\lambda$ : le faisceau $\mathrm{R}^2 \mathrm{pr}_{1!} \mathscr{L}_\lambda$ est non nul et à support dans $\{\lambda\}$. (Voir [Laumon 1987, 1.2.2.2] pour l'analogue $\ell$-adique et [Orgogozo 2003, p. 139—140] pour une discussion du même phénomène dans une autre perspective.) La flèche $\mathrm{R}^2 \mathrm{pr}_{1!} \mathscr{L}_\lambda \to \mathrm{R}^2 \mathrm{pr}_{1\star} \mathscr{L}_\lambda$ d'oubli des supports étant un isomorphisme ([Laumon 1987, 1.3.1.1]), on voit également que l'on ne dispose pas d'un analogue uniforme « naïf » du théorème [Th. finitude, 1.9 (ii)] de commutation au changement de base générique : en tout point au voisinage duquel $R^2 \operatorname{pr}_{1\star} \mathscr{L}_\lambda$ commute au changement de base, sa fibre est nécessairement nulle. Cependant, comme l'a démontré Luc Illusie — motivé par une question de Jean-Pierre Serre —, la situation est aussi bonne que possible si l'on considère des coefficients *constants* (voir [Illusie 2010]).

La formule de Grothendieck-Ogg-Šafarevič ([SGA 5 X, 7.2]) montre également que l'on ne peut obtenir de bornes uniformes sur les nombres de Betti, comme celles rappelées au paragraphe suivant, sans hypothèse de modération.

**0.4.** Dans [Katz et Laumon 1985, §3] (voir aussi par exemple [Fouvry et Katz 2001, §2] et [Katz 1990, §13]) sont établis des résultats intermédiaires entre l'énoncé topologique rappelé ci-dessus et des énoncés, à démontrer, tenant compte de la ramification sauvage. Nicholas Katz et Gérard Laumon établissent notamment le résultat suivant, sous l'hypothèse technique supplémentaire que $S$ est de type fini sur un schéma nœthérien régulier de dimension inférieure ou égale à 1 ([Katz et Laumon 1985, 3.1.2]) : il existe un entier $N \geqslant 1$ et une stratification $\mathfrak{S}$ de $S$ tels que pour chaque nombre premier $\ell$, chaque $\mathscr{K} \in \operatorname{Ob} D_c^b(X[1/N\ell], \mathbb{Q}_\ell)$ constructible le long de $\mathfrak{X}[1/N\ell]$, le complexe $Rf_\star \mathscr{K} \in \operatorname{Ob} D_c^b(S[1/N\ell], \mathbb{Q}_\ell)$ est constructible le long de $\mathfrak{S}[1/N\ell]$. (« *Constructible le long de …* » : voir 1.1.3.) De plus, les auteurs majorent uniformément la taille des images directes : il existe une constante $C$ telle que pour chaque $s \in S[1/N\ell]$ on ait l'inégalité

$$\|Rf_! \mathscr{K}\|(s) \leqslant C \cdot \sup_{x \in f^{-1}(s)} \|\mathscr{K}\|(x),$$

où $\|\mathscr{K}\|(x) \coloneqq \sum_i \dim_{\mathbb{Q}_\ell} H^i(\mathscr{K})_{\overline{x}}$ pour un choix arbitraire de point géométrique $\overline{x}$ au-dessus de $x \in X$.

Ce résultat s'appuie sur la résolution des singularités d'Hironaka et n'est intéressant que si la caractéristique générique de $S$ est nulle. Cette hypothèse permet également, grâce au lemme d'Abhyaṁkar relatif (cf. [ibid., §3.4.2] et [SGA 1 XIII, 5.5]), d'ignorer la ramification sauvage.

**0.5.** L'objectif de ce texte est de s'affranchir de ces restrictions : les « fibrations élémentaires » ([SGA 4 XI, §3]) utilisées dans [Katz et Laumon 1985] sont remplacées par des courbes nodales, dont l'ubiquité a été démontrée par Aise Johann de Jong ([de Jong 1996 ; de Jong 1997]). Armé de cet outil et familier des techniques de dévissages en cohomologie étale, la principale difficulté est de trouver un énoncé de stabilité de familles « (*uniformément*) constructibles et modérées ». On le doit à Ofer Gabber, dont la contribution à cet article ne s'arrête pas là. Un point clef est que « uniformément » doit être interprété *localement pour la topologie propre* — engendrée par la prétopologie dont les morphismes couvrants sont les morphismes propres et surjectifs — (ou la topologie $h$), idée que l'on trouve déjà dans une lettre de Richard Pink à N. Katz ([Pink 1995]), pourtant antérieure aux travaux de A. J. de Jong. (L'insuffisance d'une localisation pour la topologie finie est expliquée en 3.2.)

Les méthodes employées ici font que la plupart des résultats de cet article sont de démonstration plus simple dans le cas particulier des schémas de type fini sur un anneau de Dedekind mais nous les énonçons pour des schémas quasi-excellents, quitte à en donner la démonstration en deux temps. Il est apparu que l'hypothèse d'excellence n'était pas toujours nécessaire, notamment pour le théorème d'uniformité propre ; pour ne pas retarder à nouveau la publication de ce texte, nous avons préféré nous contenter d'esquisser les modifications à apporter (voir par exemple 3.5.5). Nous espérons que le lecteur nous pardonnera ce choix.

**0.6**.   Signalons maintenant quelques limites de la méthode. Tout d'abord, on aimerait des résultats d'*indépendance* de $\ell$, beaucoup plus forts, mais qui semblent malheureusement hors de portée : la descente cohomologique — utilisée de façon répétée — fait naturellement apparaître des différentielles dans la suite spectrale correspondante ; on ignore si elles sont « motiviques » (cf. [DE JONG 1996, introduction]). D'autre part, comme expliqué ci-dessus, on obtient notamment des majorations de nombres de Betti : peut-on calculer explicitement une borne et les nombres de Betti d'un schéma donné pour chaque $\ell$ donné ? Les techniques de cet article ne permettent pas de répondre à cette question (voir cependant [MADORE et ORGOGOZO 2015] pour le calcul des nombres de Betti). Enfin, nous ne considérons des situations sauvages que de façon superficielle : on se ramène immédiatement à des cas modérés pour lesquels les techniques classiques — voir par exemple [SGA 7 XIII, 2.1.9-11] — s'appliquent. Il serait intéressant de reprendre les résultats de cet article à la lumière des travaux plus récents sur la ramification sauvage (voir par exemple [KATO et SAITO 2013]).

**0.7**.   **Remerciements.**   C'est avec un grand plaisir que j'exprime mes sincères remerciements à Ahmed ABBES, Nick KATZ, Arthur OGUS, Joël RIOU, Jean-Pierre SERRE, ZHENG Weizhe (郑维喆) et tout particulièrement SAITO Takeshi (斎藤毅) pour avoir relevé une erreur dans une version antérieure de ce texte ainsi qu'Ofer GABBER et Luc ILLUSIE pour leur aide cruciale et leurs encouragements.

## 1. Définitions, sorites et préliminaires topologiques

**1.1**. **Constructibilité.**

**1.1.1**.   Soit $X$ un schéma nœthérien. Une **stratification** de $X$ est une partition finie (ensembliste) $\mathfrak{X}$ en sous-schémas réduits localement fermés (non vides), appelés **strates**. Si pour toute strate $S$ le bord $\partial S = \overline{S} - S$ est réunion de strates, on dit que $\mathfrak{X}$ est une **bonne stratification**. Un **schéma stratifié** est un schéma muni d'une stratification. Enfin, une stratification $\mathfrak{X}$ de $X$ est **adaptée** à une partie $Z \subseteq X$ si chaque strate de $\mathfrak{X}$ est contenue dans $Z$ ou $X - Z$.

Faisons les observations suivantes :

(1) toute stratification se raffine en une bonne stratification ([Stacks, lemme 09Y5]).

Et dans le cas d'une *bonne* stratification :

(2) toute strate, étant localement fermée, est ouverte dans son adhérence ;
(3) tout strate $S$ maximale pour la relation d'ordre partiel « $x \leqslant y$ si et seulement si $x \subseteq \overline{y}$ » est ouverte : $X \smallsetminus S = \bigsqcup_{T \neq S} T = \bigcup_{T \neq S} \overline{T}$.

Si $f : Y \to X$ est un morphisme, on note $\mathfrak{X}_Y$ ou $f^\star \mathfrak{X}$ la stratification de $Y$ déduite de $\mathfrak{X}$ par image inverse (réduite) des strates.

**1.1.2. Mise en garde.** Notons que l'image inverse d'une bonne stratification par un morphisme $f : Y \to X$ n'est pas nécessairement bonne : si $S$ est une partie localement fermée, l'inclusion $\overline{f^{-1}(S)} \subseteq f^{-1}(\overline{S})$ peut être stricte. De même, si $U$ est un ouvert d'un schéma nœthérien $X$, il n'est en général pas vrai que la stratification $\{U, X \smallsetminus U\}$ soit bonne.

Compte tenu de (1) ci-dessus, ces observations n'ont pas de conséquence fâcheuse.

**1.1.3.** Soient $X$ un schéma et $(X', \mathfrak{X}')$ un $X$-schéma stratifié nœthérien. On dit qu'un faisceau étale sur $X$ est **constructible le long de $\mathfrak{X}'$** (ou **$\mathfrak{X}'$-constructible**) si son image inverse sur $X'$ est (ensemblistement) constructible et localement constante sur chacune des strates de $\mathfrak{X}'$. Lorsque $X' = X$, on retrouve la définition de [Katz et Laumon 1985, §3.0].

Dans cet article on dira souvent « **lisse** » pour « localement constant (ensemblistement) constructible » et, assez systématiquement, « faisceau » pour « faisceau étale ».

À des fins de dévissage, le résultat suivant — obtenu en utilisant (3) *supra*[2] et procédant par récurrence — est particulièrement utile. (C'est principalement ce fait qui motive l'usage de *bonnes* stratifications.)

**1.1.4. Proposition** ([SGA 4 IX, 2.5]). *Soient $(X, \mathfrak{X})$ un schéma nœthérien stratifié et $\mathscr{F}$ un faisceau abélien constructible le long de $\mathfrak{X}$. Si $\mathfrak{X}$ est bonne, il existe une filtration finie sur $\mathscr{F}$ de gradué isomorphe à une somme directe de faisceaux du type $i_\star j_! \mathscr{L}$, où $j : S \hookrightarrow \overline{S}$ est l'immersion ouverte d'une strate de $\mathfrak{X}$ dans son adhérence, $i$ l'immersion fermée de $\overline{S}$ dans $X$ et $\mathscr{L}$ est un faisceau abélien lisse sur $S$. De plus, il existe une telle filtration de longueur majorée par un entier ne dépendant que de $\mathfrak{X}$ (et pas de $\mathscr{F}$).*

L'énoncé ci-dessus est une variante de celui de [SGA 4 IX], dans lequel on se donne une *partition* de $X$.

**1.2. Unipotence.**

**1.2.1.** Soient $X$ un schéma et $(X', \mathfrak{X}')$ un $X$-schéma stratifié nœthérien. On dit qu'un faisceau étale abélien sur $X$ est **unipotent le long de $\mathfrak{X}'$** (ou **$\mathfrak{X}'$-unipotent**) si ses tirés en arrière aux strates de $\mathfrak{X}'$ sont extensions de faisceaux *constants*. Il est dit **localement unipotent le long de $\mathfrak{X}'$** (ou **$\mathfrak{X}'$-localement unipotent**) si pour tout hensélisé strict $T$ de $X'$ — ou de façon équivalente tout $X'$-schéma strictement local $T$ —, son image inverse sur le schéma $T$ est unipotent le long de $\mathfrak{X}'_T$. Par exemple, un faisceau lisse est localement unipotent (le long de toute stratification).

(À la formulation près, cette définition se trouve dans [Pink 1995][3].)

---

[2] Comme me l'ont signalé L. Illusie et Zheng W., l'énoncé peut être mis en défaut si la stratification n'est pas bonne.

[3] Je remercie L. Illusie de m'avoir communiqué cette lettre, après la rédaction des principaux résultats de cet article.

**1.2.2.** Notons qu'un faisceau abélien constructible qui est localement unipotent le long de $\mathfrak{X}$ est automatiquement constructible le long de $\mathfrak{X}$ : une extension de faisceaux lisses est lisse ([SGA 4 IX, 2.1 (iii)]) et la lissité se teste sur les hensélisés stricts (cf. [SGA 4 IX, 2.13 (i)]).

**1.2.3.** Si $\mathscr{L}$ est un faisceau lisse sur un ouvert dense $U$ normal connexe d'un schéma $X$ et $\mathfrak{X}$ est une stratification adaptée à $U$, le faisceau constructible $j_!\mathscr{L}$ est localement unipotent le long de $\mathfrak{X}$ si et seulement si il l'est le long de la stratification $\{U, X \smallsetminus U\}$ : sur un schéma normal, l'unipotence d'un faisceau lisse se teste en les points maximaux.

**Lemme.** *Soient $X$ un schéma nœthérien normal intègre de point générique $\eta$. Un faisceau étale abélien lisse sur $X$ est unipotent, c'est-à-dire extension de faisceaux constants, si et seulement si sa restriction à $\eta$ l'est.*

**1.2.4.** Soient $G$ un groupe et $\mathscr{F}$ un $X$-faisceau étale abélien muni d'une action de $G$, c'est-à-dire d'un morphisme $G \to \mathrm{Aut}(\mathscr{F})$. Nous dirons que $\mathscr{F}$ est **$G$-unipotent** s'il existe une filtration finie croissante par des sous-faisceaux $G$-stables $\mathscr{F}_i$ telle que $G$ agisse *trivialement* sur les quotients $\mathscr{F}_{i+1}/\mathscr{F}_i$.

**1.2.5. Lemme.** *Un $G$-faisceau abélien constructible $\mathscr{F}$ sur un schéma nœthérien $X$ est $G$-unipotent si et seulement si pour tout point géométrique $\overline{x} \to X$, la fibre $\mathscr{F}_{\overline{x}}$ est $G$-unipotente.*

*Démonstration.* La condition sur les fibres géométriques est évidemment nécessaire. Considérons la réciproque. Les deux conditions étant stables par extension, restriction à un ouvert (resp. fermé) et application des endofoncteurs $j_!j^!$ et $i_\star i^\star$, un dévissage immédiat nous ramène au cas particulier où $\mathscr{F}$ est de la forme $j_!\mathscr{L}$, où $\mathscr{L}$ est un faisceau abélien lisse et $j$ une immersion ouverte. Le prolongement par zéro préservant la $G$-unipotence, on peut supposer que $\mathscr{F}$ est un faisceau lisse. En considérant la fibre en un point géométrique de $X$ (que l'on peut supposer connexe), on est ramené à montrer que si $M$ est un groupe abélien fini muni d'une action continue d'un groupe profini $\pi$ et de l'action d'un groupe $G$ commutant à l'action de $\pi$, le groupe $M$ est $G$-unipotent en tant que $\pi$-module continu si et seulement si il l'est en tant que groupe abélien. Or, si $v \in M - \{0\}$ est un point fixe de $G$, le sous-$\pi$-module $N$ de $M$ engendré par $v$ est muni d'une action triviale de $G$ et la suite exacte $0 \to N \to M \to M/N \to 0$ permet de conclure par récurrence. $\square$

## 1.3. Modération.

**1.3.1.** Un faisceau $\mathscr{F}$ sur un schéma nœthérien $X$ est dit **modéré** si pour tout hensélisé strict $T$ de $X$ — ou de façon équivalente tout $X$-schéma strictement local $T$ —, de point fermé $t$, et tout point géométrique $\overline{u}$ de $T$ localisé en $u$, l'action de $\pi_1(u, \overline{u})$ sur la fibre $\mathscr{F}_{\overline{u}}$ en $\overline{u}$ se factorise à travers un quotient d'ordre premier à l'exposant caractéristique du corps résiduel $\kappa(t)$.

**1.3.2**. Cette définition, due à O. Gabber, satisfait les deux compatibilités suivantes :

(1) un faisceau localement constant est modéré ;
(2) un faisceau sur un trait est modéré au sens ci-dessus si et seulement si sa restriction au point générique est modérée au sens usuel ([SGA 1 XIII, 2.1.1]).

**1.3.3**. Un faisceau abélien de torsion inversible sur $X$ et localement unipotent est modéré : pour tout $p$-groupe fini $P$ et tout entier $n \perp p$[①], une extension de $\mathbb{Z}/n\mathbb{Z}[P]$-modules triviaux est triviale.

**1.4. Sorites.** Soit $X$ un schéma nœthérien. La propriété, pour un faisceau abélien $\mathscr{F}$ sur $X$, d'être constructible (resp. modéré, localement unipotent) le long d'une stratification $\mathfrak{X}$ est stable par :

**1.4.1.** raffinement de la stratification : $\mathfrak{X}' = \{X'_j\}$ raffine $\mathfrak{X} = \{X_i\}$ s'il existe pour chaque $j$ un indice $i$ tel que $X'_j \subseteq X_i$ ;

**1.4.2.** image inverse du faisceau et de la stratification ;

**1.4.3.** prolongement par zéro depuis un fermé : si $\mathscr{F}$ est constructible (resp. modéré, localement unipotent) le long de $\mathfrak{X}$ et $i : X \hookrightarrow Y$ est une immersion fermée, le faisceau $i_\star \mathscr{F}$ est constructible (resp. modéré, localement unipotent) le long de la stratification de $\mathfrak{Y}$ déduite de $\mathfrak{X}$ par ajout de la partie $Y \smallsetminus X \subseteq Y$. Quitte à raffiner cette stratification $i_\star \mathfrak{X}$, il existe une *bonne* stratification de $Y$, indépendante de $\mathscr{F}$, le long de laquelle le faisceau est constructible (resp. modéré, localement unipotent).

(Signalons que la stratification ci-dessus est bonne par exemple si $X$ n'est nulle part dense dans $Y$.)

**1.4.4.** De plus, si l'on considère des faisceaux abéliens de torsion inversible sur $X$, on a également stabilité par *extension* du faisceau. C'est tautologique, sans l'hypothèse sur la torsion, pour l'unipotence locale. Cela résulte de la semi-simplicité de l'action des sous-groupes de Sylow pour la modération (cf. **1.3.3**).

**1.4.5. Remarque.** Signalons une fois pour toutes que la condition de modération est particulière en ce sens qu'elle ne fait pas intervenir de stratification. Il nous a semblé préférable, ici et par la suite, de ne pas l'isoler pour autant des deux autres conditions (constructibilité et locale unipotence, *le long d'une stratification*). Prendre garde au fait que c'est la seule des trois propriétés qui soit stable par sous-quotient quelconque ; on a cependant le substitut ci-dessous.

---

[①]On utilise dans cet article la notation de Donald Knuth pour indiquer que deux nombres sont premiers entre eux.

**1.4.6.** Soit $k : Z \hookrightarrow X$ une immersion ouverte ou fermée. Si $\mathscr{F}$ est un faisceau sur $X$ constructible (resp. localement unipotent) le long d'une stratification $\mathfrak{X}$ adaptée à $Z$ (au sens de 1.1.1), alors le sous-quotient $k_!k^\star\mathscr{F}$ du faisceau $\mathscr{F}$ est constructible (resp. localement unipotent) le long de $\mathfrak{X}$. (Amplification : dans 1.1.4, le gradué est également constructible (resp. localement unipotent) le long de $\mathfrak{X}$.)

De plus, il existe pour toute stratification $\mathfrak{X}$ une stratification la raffinant adaptée à $Z$ : intersecter les strates avec $Z$ et son complémentaire. Elle est bonne si $\mathfrak{X}$ l'est.

Signalons également le fait suivant, de vérification immédiate.

**1.4.7.** Soient $\mathfrak{X}_\alpha$ des stratifications des composantes irréductibles $X_\alpha$ de $X$. Il existe une stratification $\mathfrak{X}$ de $X$ satisfaisant la propriété suivante : tout faisceau $\mathscr{F}$ sur $X$ dont les *restrictions* aux $X_\alpha$ sont constructibles (resp. modérées, localement unipotentes) le long de $\mathfrak{X}_\alpha$ est également constructible (resp. modéré, localement unipotent) le long de $\mathfrak{X}$. (Voir 1.5.4 pour une généralisation.)

**1.4.8.** Enfin, ces conditions sont locales pour la topologie de Zariski (et même étale) : si $U_i \to X$ est un recouvrement de $X$ par des ouverts de Zariski, le faisceau $\mathscr{F}$ est constructible (resp. modéré, localement unipotent) le long de la stratification $\mathfrak{X}$ si et seulement si, pour chaque $i$, le faisceau $\mathscr{F}_{|U_i}$ l'est, le long de la stratification $\mathfrak{X}_{U_i}$ de $U_i$.

**1.4.9. Mise en garde.** Il n'est par contre pas vrai qu'un faisceau à la fois modéré sur un ouvert et sur son complémentaire est nécessairement modéré : il faut que le *prolongement par zéro* de l'ouvert à l'espace entier soit modéré. Plus précisément, si $j : U \hookrightarrow X$ est une immersion ouverte, $i : F \hookrightarrow X$ l'immersion fermée du complémentaire et $\mathscr{F}$ un faisceau de torsion inversible sur $X$ tel que $j_!j^\star\mathscr{F}$ et $i_\star i^\star\mathscr{F}$ — ou, de façon équivalente, $i^\star\mathscr{F}$ — soient modérés, alors $\mathscr{F}$ est modéré : cela résulte de la suite exacte

$$0 \to j_!j^\star\mathscr{F} \to \mathscr{F} \to i_\star i^\star\mathscr{F} \to 0.$$

(Notons que les supports des faisceaux étant disjoints, l'extension est modérée sans hypothèse sur les coefficients ; cf. 1.4.4.) Même résultat dans le cas constructible (resp. *localement* unipotent) le long d'une stratification adaptée à $U$.

Signalons cependant que cette subtilité — qui aura son importance en 4.4 — n'apparaît pas dans le cas constructible (ou *globalement* unipotent) : si $\mathscr{F}$ est un faisceau, $U$ un ouvert de complémentaire $F$ et $\mathfrak{X}$ une stratification *adaptée à* $U$, alors $\mathscr{F}$ est constructible (ou *globalement* unipotent) le long de $\mathfrak{X}$ dès que $\mathscr{F}_{|U}$ et $\mathscr{F}_{|F}$ sont constructibles (resp. globalement unipotents) le long de $\mathfrak{X}_U$ et $\mathfrak{X}_F$ respectivement.

Pour référence ultérieure, signalons la variante suivante (dont on omettra avec raison la lecture).

**1.4.10.** Amplification. Soient $j : U \hookrightarrow X$ et $k : F^\circ \hookrightarrow F := X \smallsetminus U$ deux immersions ouvertes et $\mathfrak{X}$ une stratification de $X$ adaptée à $U$ et $F^\circ$. Notons $j'$ l'immersion ouverte de $U' := U \cup F^\circ$ dans $X$. Alors, si $j_!j^\star\mathscr{F}$ est constructible (resp. modéré, localement unipotent) le long de $\mathfrak{X}$ et $k_!k^\star(\mathscr{F}_{|F})$ est constructible (resp. modéré, localement unipotent) le long de $\mathfrak{X}_F$, alors le faisceau $j'_! j'^\star \mathscr{F}$ est également constructible (resp. modéré, localement unipotent) le long de $\mathfrak{X}$.

**1.5. Constructibilité, modération et unipotence uniformes.** Soient $X$ un schéma nœthérien, $\mathsf{T}$ une topologie sur la catégorie des $X$-schémas de type fini et $\mathbb{F} = (\mathscr{F}_\lambda)_{\lambda \in \Lambda}$ une famille de faisceaux étales abéliens sur $X$ (ou plus généralement sur des ouverts $U_\lambda$ de $X$ ; voir 1.5.2).

**1.5.1.** La *famille* $\mathbb{F}$ est dite **constructible** (resp. **modérée**, **localement unipotente**) **pour la topologie** $\mathsf{T}$ s'il existe une famille couvrante $(q_i : X_i \to X)_{i \in I}$ pour la topologie $\mathsf{T}$ et des stratifications $\mathfrak{X}_i$ des $X_i$ telles que les faisceaux $q_i^\star \mathscr{F}_\lambda$, $\lambda \in \Lambda$, soient tous constructibles (resp. modérés, localement unipotents) le long des $\mathfrak{X}_i$ (1.1.3). Si l'on souhaite être plus précis, on dit également les faisceaux $\mathscr{F}_\lambda$ sont **modérés** (resp. **rendus localement unipotents**) **par les** $q_i$, $\mathfrak{X}_i$ (ou simplement par les $q_i$). Pour insister sur le fait que l'on considère une *famille* de faisceaux et non un seul faisceau, on dira parfois que le famille est « *uniformément* » constructible (resp. modérée, localement unipotente).

Lorsque $\mathsf{T}$ est la topologie grossière, on dit simplement qu'une famille est « constructible (resp. modérée, localement unipotente) *sur $X$* ».

Signalons que la définition précédente, et l'introduction du terme de « famille » de faisceaux, n'est qu'une commodité de langage dont on pourrait se passer (comparer avec [Katz et Laumon 1985, §3]).

**1.5.2.** *Convention.* Pour tout entier $n$, on note « $- \mapsto -[1/n]$ » le foncteur de changement de base de $\mathrm{Spec}(\mathbb{Z})$ à $\mathrm{Spec}(\mathbb{Z}[1/n])$. Par exemple, si $X$ est un schéma, on note $X[1/n]$ le produit fibré $X \times_{\mathrm{Spec}(\mathbb{Z})} \mathrm{Spec}(\mathbb{Z}[1/n])$. De même pour un morphisme de schémas ou une stratification par exemple. Si cela ne semble pas prêter à confusion, nous noterons parfois « $- \mapsto -^\dagger$ » ce foncteur, l'entier $n$ étant alors sous-entendu.

Pour éviter certaines lourdeurs, nous omettrons parfois la mention de ces foncteurs. Par exemple, étant donné un schéma nœthérien $X$, on dira de manière abusive qu'une famille de faisceaux $\mathscr{F}_\lambda$ de $\mathbb{Z}/n_\lambda\mathbb{Z}$-modules sur $X[1/n_\lambda]$ est « constructible sur $X$ » s'il existe une stratification $\mathfrak{X}$ de $X$ telle que les $\mathscr{F}_\lambda$ soient constructibles le long de $\mathfrak{X}[1/n_\lambda]$.

**1.5.3.** Une famille est constructible (resp. modérée, localement unipotente) pour la topologie étale si et seulement si elle l'est pour la topologie grossière (c'est-à-dire sur $X$).

La proposition suivante — dont nous ne ferons pas usage dans cet article — limite sérieusement l'intérêt de la définition 1.5.1 dans le cas constructible. (Voir également 1.6.7 *infra* pour un énoncé de modération « automatique » pour une famille à un élément.) Elle renforce le résultat bien connu selon lequel un faisceau $\mathscr{F}$ sur $X$ est constructible si et seulement si son image inverse par un morphisme surjectif de présentation finie l'est ([SGA 4 IX, 2.8] ; voir aussi [Orgogozo 2006, §9.1, 10.5]).

**1.5.4**. **Proposition.** *Soient $X$ un schéma nœthérien et $(X', \mathfrak{X}')$ un $X$-schéma surjectif de type fini, stratifié. Il existe une stratification $\mathfrak{X}$ de $X$ telle que tout faisceau $\mathscr{F}$ sur $X$ constructible le long de $\mathfrak{X}'$ soit constructible le long de $\mathfrak{X}$.*

Une famille constructible localement pour les topologies considérées dans cette article (topologie $h$, topologie des altérations) est donc même constructible pour la topologie grossière.

*Démonstration.* Supposons $X$ irréductible pour simplifier (cf 1.4.7). Il existe un ouvert $V \subseteq X'$ d'une strate de $\mathfrak{X}'$ tel que $U \coloneqq f(V)$ contienne le point générique de $X$. Quitte à rétrécir $V$, on peut supposer $U$ ouvert et le faisceau constructible $\Pi \coloneqq (V \to U)_\star \underline{2}$ — où $\underline{2}$ est le faisceau constant de valeur $\{0, 1\}$ — *lisse*, de formation compatible au changement de base.

Montrons que sous ces hypothèses tout faisceau $\mathscr{F}$ sur $X$ lisse sur $V \subseteq X'$ est également lisse sur $U \subseteq X$ ; ceci permet de conclure par récurrence nœthérienne. Le schéma $U$ étant nœthérien donc localement connexe par arcs ([SGA 4 IX, 2.12]), il suffit de vérifier que pour chaque morphisme $\overline{y} \rightsquigarrow \overline{x}$ de points géométriques de $U$, le morphisme de spécialisation $\mathscr{F}_{\overline{x}} \to \mathscr{F}_{\overline{y}}$ est un isomorphisme ([SGA 4 IX, 2.13]). Par hypothèse sur le faisceau $\Pi$, les flèches $\pi_0(V_{\overline{y}}) \to \pi_0(V \times_U U_{(\overline{x})}) \leftarrow \pi_0(V_{\overline{x}})$ sont des isomorphismes de sorte que toute composante connexe $W$ de $V \times_U U_{(\overline{x})}$ rencontre $V_{\overline{y}}$ et $V_{\overline{x}}$. Ainsi, le morphisme $\overline{y} \rightsquigarrow \overline{x}$ se relève dans $W$ en une chaîne de spécialisations-générisations ; la conclusion résulte alors de la lissité de $\mathscr{F}$ sur $W$. □

## 1.6. Un « lemme » d'Ofer Gabber.

**1.6.1**. *Motivations.* Comme signalé dans l'introduction, nous utiliserons principalement dans cet article la *topologie propre* ou — compte tenu du caractère Zariski (et même étale) local de nos conditions — la *topologie $h$*. (Voir aussi 1.6.5 *infra*.) Les dévissages, de nature géométrique ou faisceautique (cf. 1.1.4), et le formalisme des six opérations font jouer un rôle clef aux immersions fermées. Le résultat suivant de O. Gabber a pour conséquence que tout se passe dans cet article comme si *tout $h$-recouvrement d'un fermé est trace d'un $h$-recouvrement de l'espace ambiant.*

**1.6.2**. **Proposition** (O. Gabber). *Soient $X$ un schéma nœthérien, $i : F \hookrightarrow X$ une immersion fermée et $f : F' \to F$ un morphisme fini (resp. propre) et surjectif. Il existe un morphisme fini surjectif (resp. une altération) $X' \to X$ tel que toute composante irréductible $T$ de $F \times_X X'$, munie de la structure réduite, soit la source d'un $F$-morphisme vers $F'$.*

$$
\begin{array}{ccc}
F' & \cdots\cdots & \\
\downarrow & & \\
F & \llleftarrow\cdots F_{X'} & \llleftarrow \coprod_\alpha T_\alpha \\
i \downarrow & \square & \uparrow \\
X & \llleftarrow\cdots X' &
\end{array}
$$

La suite de cette section 1.6 est consacrée à la démonstration de cette proposition et quelques conséquences.

**1.6.3**. *Cas propre.* Il résulte de [Orgogozo 2006, lemme 4.3] (également dû à O. Gabber) que le cas propre se ramène au cas fini, traité ci-dessous. Pour la commodité du lecteur, nous rappelons brièvement l'argument. Comme expliqué dans [ibid.], il résulte du théorème de platification par éclatement, qu'il existe pour tout sous-schéma fermé intègre $Z$ de $F$ un éclatement $\widetilde{Z} \to Z$ et un diagramme commutatif

$$
\begin{array}{ccccc}
F' & \dashleftarrow\!\!\!\!-\!\!\!\dashrightarrow & Z' & \dashleftarrow & \widetilde{Z'} \\
{\scriptstyle\text{propre}}\downarrow & & {\scriptstyle\text{altération}}\downarrow & & \downarrow{\scriptstyle\text{fini}} \\
F & \longleftarrow\!\!\!\!-\!\!\!\!\longrightarrow & Z & \dashleftarrow & \widetilde{Z} = \text{Écl}_R(Z) \,.
\end{array}
$$

Par compacité de $F^{\text{cons}}$, il existe un nombre fini de fermés intègres $Z_i$ tels que $F$ soit la réunion des ouverts $Z_i - R_i$, où $R_i$ est le centre de l'éclatement $\widetilde{Z_i} \to Z_i$.

Soit $\widetilde{X}$ l'éclatement dominant les éclatés de $X$ le long des $R_i$ et notons $T_\alpha$ les composantes irréductibles du produit fibré $F_{\widetilde{X}} \coloneqq F \times_X \widetilde{X}$. Par construction, chaque $T_\alpha$ s'envoie sur un $\widetilde{Z_i}$. Le schéma $T' = \amalg_\alpha T'_\alpha$, obtenu par changement de base comme dans le diagramme ci-dessous, est alors *fini surjectif* sur $F_{\widetilde{X}}$ et il existe un $F$-morphisme de $T'$ vers $F'$. Comme annoncé, ceci nous ramène au cas où le morphisme $F' \twoheadrightarrow F$ est *fini* (surjectif).

$$
\begin{array}{ccccc}
F' & \longleftarrow & \widetilde{Z_i'} & \dashleftarrow & T'_\alpha \\
{\scriptstyle\text{propre}}\downarrow & & \downarrow{\scriptstyle\text{fini}} & \square & \downarrow{\scriptstyle\text{fini}} \\
F & \longleftarrow & \widetilde{Z_i} & \dashleftarrow & T_\alpha \subseteq F \times_X \widetilde{X} \\
\downarrow & & \downarrow & & \uparrow \\
& & \text{Écl}_{R_i}(X) & & \\
\downarrow & & & & \uparrow \\
X & \longleftarrow & & & \widetilde{X} = \text{Écl}_{\bigcup_i R_i}(X)
\end{array}
$$

**1.6.4**. *Cas fini.* On commence par supposer $X$ affine pour simplifier.

**1.6.4.1**. Rappelons ([Bourbaki A, VI, §6, n°5], [Lombardi et Quitté 2011, III.§4, VII.§4]; voir aussi [Bhargava et Satriano 2014, §6]) que si $A$ est un anneau et $f \in A[t]$ est un polynôme unitaire, on peut définir l'**algèbre de décomposition universelle** $\text{Adu}_A(f)$ de $f$, quotient de l'anneau $A[y_i, 1 \leqslant i \leqslant \deg(f)]$ par les relations obtenues en imposant l'égalité $\prod_i (t - y_i) = f(t)$. Si $P$ est une famille finie de polynôme unitaires, on note $\text{Adu}_A(P)$ le produit tensoriel des $\text{Adu}_A(f)$, pour $f \in P$, au-dessus de $A$. Lorsque $X = \text{Spec}(A)$, on note $\text{Sdu}_X(P)$ le **schéma de décomposition universel**, spectre de cette algèbre. Par construction, pour tout point $x'$ de $\text{Sdu}_X(P)$ d'image $x$ dans $X$, l'extension $\kappa(x')/\kappa(x)$ est finie quasi-galoisienne (=normale).

**1.6.4.2**. Lorsque dans **1.6.2** le schéma $X$ est affine, le morphisme fini $F' \to F$ correspond à un morphisme fini $\mathcal{O}(F) \to \mathcal{O}(F')$ de $\mathcal{O}(X)$-algèbres. Choisissons arbitrairement des générateurs $a_1, \ldots, a_n$ de $\mathcal{O}(F')$ sur $\mathcal{O}(F)$ satisfaisant des équations unitaires $q_i \in Q \subseteq \mathcal{O}(F)[t]$, que l'on relève (arbitrairement) en des équations unitaires $p_i \in P \subseteq \mathcal{O}(X)[t]$. Soit $X' \coloneqq \mathrm{Sdu}_X(P)$ le schéma de décomposition universel associé ; le produit fibré $F_{X'} \coloneqq F \times_X X'$ est alors $\mathrm{Sdu}_F(Q)$.

Soit $x'$ un point de $X'$ d'image $x \in F$ dans $X$. Pour tout point $f'$ de $F'$ au-dessus de $x$, il existe par construction un $x$-morphisme $x' \to f'$.

$$\begin{array}{ccccc}
x' & & F_{X'} & \hookrightarrow & X' \\
\vdots & & & & \\
\downarrow & & & & \\
f' & & F' & \square & \\
& & & & \\
\downarrow & & \downarrow & & \downarrow \\
x & & F & \hookrightarrow & X
\end{array}$$

Appliquant ceci aux points maximaux de $F_{X'}$, on observe que pour toute composante irréductible $T_\alpha$ de $F_{X'}$ il existe un $F$-morphisme de son point générique vers $F'$. Par construction également, ce morphisme s'étend, de façon *unique*, à la composante irréductible $T_\alpha$ tout entière (munie de la structure réduite). En effet, tout diagramme commutatif comme-ci dessous, où $K, K'$ et $K''$ sont les corps des fractions d'anneaux intègres $R, R'$ et $R''$, se complète par un $R$-morphisme $R' \to R''$.

$$\begin{array}{ccccc}
K'' & \longleftarrow & R'' & \twoheadleftarrow & \mathrm{Adu}_R(q_1, \ldots, q_n) \\
\uparrow & & \uparrow & & \\
K' & \longleftarrow & R' & = & R[x_1, \ldots, x_n],\ q_i(x_i) = 0 \\
\uparrow & & \uparrow & & \\
K & \longleftarrow & R & &
\end{array}$$

**1.6.4.3**. Cas général : on ne suppose plus $X$ affine, mais on se donne un recouvrement fini par des ouverts affines $U_i$. Posons $F_i \coloneqq F \times_X U_i$.

Lorsque $F_i$ est non vide, on considère un schéma de décomposition universel $U'_i \twoheadrightarrow U_i$ associé (non canoniquement) comme ci-dessus au diagramme $F'_i \twoheadrightarrow F_i \hookrightarrow U_i$ (où $F'_i \coloneqq F' \times_F F_i$). Le morphisme $U'_i \to U_i$ s'étend en un morphisme fini $X'_i \twoheadrightarrow X$ : considérer le coproduit de $X \smallsetminus U_i$ avec une approximation de la normalisation de $X$ dans $U'_i$, qui coïncide avec $U'_i$ sur $U_i$. (Écrire la normalisation comme une limite de $X$-schémas finis.)

Lorsque $F_i = \emptyset$ (c'est-à-dire $F \cap U_i = \emptyset$), on pose $X'_i \coloneqq X$.

Soit $X'$ le produit fibré sur $X$ des $X'_i$. La démonstration du paragraphe précédent s'applique *mutatis mutandis*, l'unicité susmentionnée des prolongements nous permettant de se restreindre aux schémas affines $U_i$.

**1.6.5.** La proposition 1.6.2, dans le cas particulier où $F = X$, entraîne qu'une famille constructible et modérée (resp. constructible et localement unipotente) pour la topologie $h$ l'est aussi pour la topologie des altérations.

**1.6.6. Remarque** (O. Gabber). Bien que nous n'en ferons pas usage, observons que si une famille est « modérée » (ou rendue localement unipotente) par une altération, elle l'est aussi par une altération *génériquement étale*. Cela résulte de l'invariance du topos étale par homéomorphisme universel ([SGA 4 VIII, §1]) et du fait que toute altération $X' \to X$ se décompose en $X' \to X'' \to X$, où $X'' \to X$ est une altération génériquement étale et $X' \to X''$ est un morphisme fini radiciel. En effet, si $X' \to X$ est localement défini par une injection d'anneaux $A \hookrightarrow B$, avec $A$ et $B$ intègres, on peut factoriser cette injection en $A \hookrightarrow A[B^{p^n}] \hookrightarrow B$ ; cette construction se recolle et, pour $n$ grand, l'extension $\mathrm{Frac}(A[B^{p^n}]) / \mathrm{Frac}(A)$ ainsi obtenue est étale (cf. p. ex. [Bourbaki A, V §7 n°7, cor. 2 (démonstration)]).

Terminons par une application, qui est une conséquence immédiate du lemme de Gabber précédent et de la possibilité de trivialiser un faisceau lisse par un morphisme fini étale.

**1.6.7. Proposition.** *Tout faisceau constructible sur un schéma nœthérien est modéré par un morphisme* fini *surjectif.*

*Démonstration.* Quitte à plonger le faisceau, $\mathscr{F}$, dans $\mathbb{Z}/n\mathbb{Z}[\mathscr{F}]$ pour un entier $n > 1$, on peut supposer le faisceau abélien. (On utilise le fait que $\mathscr{F} \mapsto \mathbb{Z}/n\mathbb{Z}[\mathscr{F}]$ commute aux fibres.) Par passage à la limite — utiliser par exemple [Thomason et Trobaugh 1990, C.9] et [SGA 4 IX, 2.7.4] — et stabilité par image inverse, on peut supposer le schéma nœthérien *excellent*. Par dévissage (triangle $j_! j^\star \to \mathrm{Id} \to i_\star i^\star \xrightarrow{+1}$) et la proposition 1.6.2, on se ramène par récurrence nœthérienne sur $X$ au cas d'un faisceau de la forme $j_! \mathscr{L}$, où $j : U \hookrightarrow X$ est une immersion ouverte d'image dense et $\mathscr{L}$ est un faisceau lisse. Soit $U' \to U$ un revêtement étale trivialisant $\mathscr{L}$ et considérons le normalisé $X' \to X$ de $X$ dans $U'$. Le faisceau $(X' \to X)^\star (j_! \mathscr{L})$ est modéré. $\square$

## 1.7. Extensions ponctuelles.

**1.7.1.** Soient $f : X \to Y$ un morphisme de topos, $\{\mathscr{F}_\lambda\}_{\lambda \in \Lambda}$ des faisceaux abéliens sur $X$ et $N \geqslant 0$ un entier. On dit qu'un faisceau $\mathscr{G}$ sur $Y$ est $N$-**extension ponctuelle de sous-quotients des** $\mathscr{F}_\lambda$ (resp. **relativement à** $f$) s'il existe pour chaque point $y$ de $Y$ une filtration à $N$ crans de la fibre $\mathscr{G}_y$ dont les quotients successifs sont des sous-quotients de fibres des $\mathscr{F}_\lambda$ en des points de $X$ (resp. au-dessus de $y$).

**1.7.2.** *Sorites.* Les faits suivants sont de vérification immédiate.
   (i) Soient $f : X \to Y$ un morphisme, $\mathscr{F}$ un faisceau abélien sur $X$, un morphisme $\pi : X' \to X$ et $f' : X' \to Y$ le composé $f \circ \pi$. Si un faisceau $\mathscr{G}$ sur $Y$ est $N$-extension ponctuelle de sous-quotients de $\mathscr{F}' \coloneqq \pi^\star \mathscr{F}$ relativement à $f'$, il est également $N$-extension ponctuelle de sous-quotients de $\mathscr{F}$ relativement à $f$.

(ii) Soient $f : X \to Y$ un morphisme, $k : Z \hookrightarrow X$ une immersion, et $\mathscr{F}$ un faisceau abélien sur $Z$. Un faisceau $\mathscr{G}$ sur $Y$ est $N$-extension ponctuelle de sous-quotients de $\mathscr{F}$ relativement au composé $g \coloneqq f \circ k$ si et seulement si il est $N$-extension ponctuelle de $k_!\mathscr{F}$ relativement à $f$.

(iii) Soient $f : X \to Y$ un morphisme, $\mathscr{F}$ un faisceau sur $X$,

$$\mathscr{K}_1 \to \mathscr{K} \to \mathscr{K}_2 \xrightarrow{+1}$$

un triangle distingué de complexes sur $Y$ et $q$ un entier. Si chaque $\mathrm{H}^q(\mathscr{K}_\lambda)$, $\lambda \in \{1, 2\}$, est $N_\lambda$-extension ponctuelle de sous-quotients de $\mathscr{F}$ relativement à $f$, le faisceau $\mathrm{H}^q\mathscr{K}$ est $(N_1 + N_2)$-extension ponctuelle de sous-quotients de $\mathscr{F}$ relativement à $f$.

(On appliquera ceci au cas particulier où les deux complexes $\mathscr{K}_\lambda$ sont obtenus par application du foncteur $\mathrm{R}f_\star$ à une suite exacte $0 \to \mathscr{F}_1 \to \mathscr{F} \to \mathscr{F}_2 \to 0$ de faisceaux sur $X$.)

(iv) Soient $f : X \to Z = (g : Y \to Z) \circ (h : X \to Y)$ un morphisme, $\mathscr{F}$ un faisceau sur $X$, $\mathscr{G}$ un faisceau sur $Y$ et $\mathscr{H}$ un faisceau sur $Z$. Si $\mathscr{G}$ est $N_1$-extension ponctuelle de sous-quotients de $\mathscr{F}$ (relativement à $h$) et $\mathscr{H}$ est $N_2$-extension ponctuelle de sous-quotients de $\mathscr{G}$ (relativement à $g$), alors $\mathscr{H}$ est $(N_1 \times N_2)$-extension ponctuelle de sous-quotients de $\mathscr{F}$ (relativement à $f$).

(On appliquera ceci au cas particulier où $\mathscr{G} \coloneqq \mathrm{R}^i h_\star \mathscr{F}$ et $\mathscr{H} \coloneqq \mathrm{R}^j g_\star \mathscr{G}$.)

(v) Soient

$$\begin{array}{ccc} X & \xleftarrow{g'} & X' \\ f \downarrow & & \downarrow f' \\ Y & \xleftarrow{g} & Y' \end{array}$$

un carré commutatif avec $g$ *surjectif* sur les points et $\mathscr{F}$ (resp. $\mathscr{G}$) un faisceau abélien sur $X$ (resp. $Y$). Pour tout entier $N \geqslant 0$, le faisceau $\mathscr{G}$ est $N$-extension ponctuelle de sous-quotients de $\mathscr{F}$ relativement à $f$ si et seulement si le faisceau $g^\star \mathscr{G}$ est $N$-extension ponctuelle de sous-quotients de $g'^\star \mathscr{F}'$ relativement à $f'$.

(On appliquera ceci au cas particulier d'un carré cartésien de morphismes munis de la topologie étale, $f$ est un morphisme *propre*, $\mathscr{G} \coloneqq \mathrm{R}^q f_\star \mathscr{F}$, et $g^\star \mathscr{G} = \mathrm{R}^q g'_\star \mathscr{F}'$.)

## 2. Lemmes d'Abhyankar et courbes nodales

### 2.1. Paires toriques et courbes nodales : rappels.

**2.1.1.** Rappelons que si $X$ est un schéma et $j : U \hookrightarrow X$ est un ouvert dense, la paire $(X, U)$ est dite **torique** ([Mochizuki 1999, 1.2]) ou **log-régulière** ([STG VI, 1.4]) si le schéma $X$ muni de la log-structure[5] image directe $\mathscr{O}_X \cap j_\star \mathscr{O}_U^\times \hookrightarrow \mathscr{O}_X$ est **log-régulier** et $U$ est l'ouvert de trivialité de cette log-structure. (Par définition, un log-schéma log-régulier est fin et saturé [fs].) Tout log-schéma $(X, M_X)$ log-régulier est de ce type ; cf. p. ex. [GABBER et RAMERO 2016, 10.5.53]. L'exemple typique d'une paire torique est constitué d'un schéma régulier $X$ et du complémentaire $U$ du support d'un diviseur à croisements normaux (resp. à croisements normaux strict). Une telle paire est dite **régulière** (resp. **strictement régulière**). Réciproquement, toute paire torique $(X, U)$ avec $X$ *régulier* est une paire régulière au sens précédent. Nous commettrons l'abus de notation suivant : si $(X, U)$ est torique, on désignera également par $(X, U)$ le *log-schéma* muni de la log-structure introduite ci-dessus.

---

[5] Comme expliqué en [Gabber et Ramero 2016, 7.2.12] cette log-structure sur $X$, muni de la topologie étale, n'est pas nécessairement fine.

**2.1.2.** Soit maintenant $X$ un schéma nœthérien. Rappelons qu'un $X$-schéma $Y$ est une **courbe nodale** si le morphisme structural $h : Y \to X$ est projectif, plat, à fibres géométriques des courbes connexes ayant au pire des singularités quadratiques ordinaires. Nous dirons d'une courbe nodale $h : Y \to X$ qu'elle est **adaptée** à une paire $(Y°, X°)$, où $Y°$ (resp. $X°$) est un ouvert dense de $Y$ (resp. $X$) lorsque les conditions suivantes sont satisfaites :

— le morphisme $h$ est lisse au-dessus de $X°$ ;
— il existe un diviseur $D$ étale sur $X$ contenu dans le lieu lisse de $f$ tel que $Y° = h^{-1}(X°) \cap (Y - D)$.

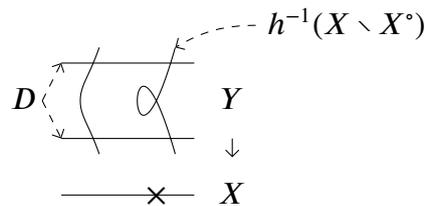

Un lien avec la géométrie logarithmique est donné par la proposition bien connue suivante ([STG VI, 1.9], [MOCHIZUKI 1999, 4.2] ou [SAITO 2004, 1.12]).

**2.1.3. Proposition.** *Soit* $h : Y \to X$ *une courbe nodale adaptée à une paire d'ouverts denses* $(Y°, X°)$ *telle que* $(X, X°)$ *soit* torique. *Alors le morphisme* $(Y, Y°) \to (X, X°)$ *est* log-lisse *et la paire* $(Y, Y°)$ *est torique.*

(Voir aussi 2.3.6.)

## 2.2. Pureté.

**2.2.1.** Soient $X$ un schéma régulier connexe, $D$ un diviseur à croisements normaux et $U = X \smallsetminus D$. Pour tout faisceau $\mathscr{L}$ lisse sur $U = X \smallsetminus D$, il est équivalent de demander que le prolongement par zéro $j_!\mathscr{L}$ soit modéré au sens de 1.3.1 ou bien que $\mathscr{L}$ soit *modéré le long de $D$* au sens usuel (c'est-à-dire : pour tout point maximal $d$ de $D$, l'inertie sauvage du corps discrètement valué $\mathrm{Frac}(\mathscr{O}_{X,d})$ agit trivialement). Seule une implication est non triviale ; elle résulte immédiatement du lemme d'Abhyaṅkar usuel ([SGA 1 XIII, 5.2]).

**2.2.2.** Considérons maintenant une variante logarithmique du critère de modération précédent s'appuyant sur un théorème de pureté logarithmique (Katô K., cf. p. ex. [Mochizuki 1999, 3.3]). On ne suppose plus la paire $(X, U)$ régulière mais seulement *torique*. La condition de modération en les points maximaux du diviseur $D$ garde un sens, car $X$ est normal ([Kato 1994, 4.1] ou [Gabber et Ramero 2016, 10.5.29]). Si un faisceau abélien $\mathscr{L}$ de $\mathbb{Z}/n\mathbb{Z}$-modules, lisse sur $U$, satisfait cette condition, il existe un revêtement étale $U' \to U$, modéré le long de $D$ (au sens usuel, [ibid., 11.3.16], [Illusie 1981, 2.6]) trivialisant $\mathscr{L}$. Soit $X'$ le normalisé de $X$ dans $U'$; le morphisme $X' \to X$ est fini et surjectif. D'après [Gabber et Ramero 2016, 11.3.43], la paire $(X', U')$ est torique et le morphisme de log-schémas $(X', U') \to (X, U)$ est un revêtement log-étale, au sens de [ibid., 11.3.13]. Lorsque $X$ est strictement local d'exposant caractéristique résiduel $p$, le groupe fondamental logarithmique de la paire $(X, U)$ est abélien d'ordre profini premier à $p$ ([ibid., 11.3.41]). Il en résulte formellement que le faisceau $j_!\mathscr{L}$ est modéré. Résumons ces observations sous la forme d'un lemme.

**2.2.3. Lemme.** *Soient $(X, U)$ une paire* torique, *$j : U \hookrightarrow X$ l'immersion ouverte et $\mathscr{L}$ un faisceau abélien lisse sur $U$ de $\mathbb{Z}/n\mathbb{Z}$-modules ($n \geq 1$). Les conditions suivantes sont équivalentes :*

  (i) *le prolongement par zéro $j_!\mathscr{L}$ est modéré (au sens de 1.3.1);*
  (ii) *$\mathscr{L}$ est modéré le long de $D$ : pour tout point maximal $d$ de $D$, l'inertie sauvage du corps discrètement valué $\mathrm{Frac}(\mathcal{O}_{X,d})$ agit trivialement sur la fibre de $\mathscr{L}$ en un point géométrique au-dessus de $d$.*
  (iii) *il existe un revêtement log-étale $(X', U') \to (X, U)$ tel que $\mathscr{L}_{|U'}$ soit trivial.*

**2.2.4. Remarque.** L'analogue de ce lemme pour l'unipotence est faux : il existe par exemple une paire torique $(X, U)$ avec $X$ strictement local de dimension 2 et un faisceau $\mathscr{L}$ lisse sur $U$ qui n'est *pas* unipotent (sur $U$) mais localement unipotent en les points maximaux de $X - U$. En effet, comme expliqué en [ibid., 11.3.52-55], on peut calculer le groupe fondamental du spectre épointé $X^\star := X - \{x\}$ (pour $X$ de dimension quelconque) et en particulier ([ibid., 11.3.57]) construire une paire comme ci-dessus telle que $\pi_1(X^\star) = \mathbb{Z}/2\mathbb{Z}$.

**2.3. Image directe : courbe nodale.** Faute de référence adéquate, ce paragraphe est consacré à la démonstration de la proposition suivante. (Signalons cependant que, dans le cas localement unipotent, cet énoncé et une esquisse de preuve — différente — figurent dans [Pink 1995].)

**2.3.1. Proposition.** *Soient $X$ un schéma nœthérien, $h : Y \to X$ une courbe nodale adaptée à une paire d'ouverts denses $(Y^\circ, X^\circ)$, avec $X^\circ$ normal, et $\mathscr{L}$ un faisceau lisse de $\mathbb{Z}/n\mathbb{Z}$-modules sur $Y^\circ$, où $n$ est un entier inversible sur $X$, tel que $j_!\mathscr{L}$ soit modéré (resp. localement unipotent le long de la stratification $\mathfrak{Y} := \{Y^\circ, Y \smallsetminus Y^\circ\}$ de $Y$), où $j$ désigne l'immersion ouverte $Y^\circ \hookrightarrow Y$. Alors, les faisceaux $\mathrm{R}^i h_\star(j_!\mathscr{L})$ ($0 \leq i \leq 2$) sont constructibles et modérés (resp. localement unipotents) le long de la stratification $\mathfrak{X} := \{X^\circ, X \smallsetminus X^\circ\}$ de $X$.*

Nous traitons d'abord le cas où la paire $(X, X^\circ)$ est torique : cela est suffisant pour établir le théorème de constructibilité uniforme pour les morphismes propres dans le cas des schémas quasi-excellents. Le cas des schémas nœthériens généraux est ensuite esquissé.

Fixons un entier $i \in \{0, 1, 2\}$.

**2.3.2.** *Constructibilité*. Le faisceau $R^i h_\star(j_! \mathscr{L})$ est nul sur le fermé $X \smallsetminus X^\circ$ : cela résulte du théorème de changement de base propre et du fait que l'ouvert $Y^\circ$ ne rencontre pas le fermé $h^{-1}(X \smallsetminus X^\circ)$. Montrons maintenant que $R^i h_\star(j_! \mathscr{L})$ est lisse sur l'ouvert $X^\circ$. On peut supposer $X = X^\circ$. Le lemme bien connu ci-dessous et la propreté de $h$ nous ramènent au cas où $X$ est un *trait* hensélien, qui résulte de [SGA 7 XIII, 2.1.11].

(Le même argument se trouve dans [Pink 1995, théorème 4].)

**Lemme.** *Soient $X$ un schéma nœthérien, $(f : Y \to X, \mathscr{F})$ une paire cohomologiquement propre et $c \geqslant 0$ un entier tel que $\mathscr{H} \coloneqq R^c f_\star \mathscr{F}$ soit constructible. Alors, $\mathscr{H}$ est lisse si et seulement si pour tout $X$-trait strictement hensélien $T$, le faisceau $R^c f_{T\star} \mathscr{F}_{|Y_T}$ est lisse sur $T$.*

*Démonstration.* Par propreté cohomologique, il s'agit de montrer que le faisceau constructible $\mathscr{H}$ est lisse si ses tirés en arrière sur tout $X$-trait strictement hensélien le sont. Or, d'après [SGA 4 IX, 2.13 (i)]), il faut montrer que pour toute spécialisation $b \rightsquigarrow a$ de points géométriques de $X$, le morphisme $\mathscr{H}_a \to \mathscr{H}_b$ est un isomorphisme. D'après [ÉGA II, 7.1.9], il existe un morphisme local $T \to X_{(a)}$, induisant un isomorphisme sur les corps résiduels en les points fermés et envoyant le point générique de $T$ sur (l'image de) $b$. Ceci suffit pour conclure. $\square$

**2.3.3.** *Modération : cas où la paire $(X, X^\circ)$ est torique.* Soit $\pi : (Y', Y'^\circ) \to (Y, Y^\circ)$ un revêtement log-étale tel que $\mathscr{L}_{|Y'^\circ} \coloneqq \pi^{\circ\star} \mathscr{L}$ soit constant (2.2.3). La transitivité des images directes et un dévissage classique (cf. par exemple [Th. finitude, A, 1.3.3]), reposant sur le fait que le prolongement par zéro du conoyau de la flèche $\mathscr{L} \hookrightarrow \pi^\circ_\star \pi^{\circ\star} \mathscr{L}$ est également modéré, nous ramènent à montrer que pour tout revêtement log-étale $\tau : (Y', Y'^\circ) \to (Y, Y^\circ)$, le faisceau $R^i(h \circ \tau)_\star (j'_! \mathbb{Z}/n\mathbb{Z})$ est modéré. (On note $j'$ le changement de base $j \times_Y Y'$.) D'après les résultats rappelés au paragraphe précédent, ce faisceau est le prolongement par zéro d'un faisceau *lisse* sur $X^\circ$. Le lemme d'Abhyaṁkar logarithmique (2.2.3), nous permet de supposer dans le cas torique que $X$ est un *trait* strictement hensélien, dont nous notons $\eta$ le point générique de sorte que $X^\circ = \{\eta\}$, sans quoi $X = X^\circ$ et il n'y a rien à démontrer. Soit $\overline{\eta}$ un point géométrique au-dessus de $\eta$ ; on veut montrer que l'action de $\pi_1(\eta, \overline{\eta})$ sur $H^i_c(Y'^\circ_{\overline{\eta}}, \mathbb{Z}/n\mathbb{Z})$ est modérée. Par lissité de $Y'^\circ_\eta$ sur $\eta$ et dualité de Poincaré, il suffit de démontrer la modération de la cohomologie sans support (en degré cohomologique $2-i$). Ceci est un cas particulier du résultat suivant ([Illusie 2002, 8.4.3]) :

> Soient $(W, \mathscr{M}_W)$ un log-schéma fin et saturé, propre et log-lisse sur $(X, \eta)$ d'ouvert de trivialité $W^\circ$ contenu dans $W_\eta$, et $n$ un entier inversible sur $X$. Alors l'action de $\pi_1(\eta, \overline{\eta})$ sur chaque $H^j(W^\circ_{\overline{\eta}}, \mathbb{Z}/n\mathbb{Z})$ est modérée.

**2.3.4**. *Unipotence : cas d'un trait.* Commençons par le cas où $X$ est un trait strictement hensélien et $X° = \{\eta\}$. Notons $s$ le point fermé de $X$ et $\overline{\eta}$ un point géométrique au-dessus de $\eta$. On veut montrer que l'action de $\pi_1(\eta, \overline{\eta})$ sur $\mathrm{H}^i_c(Y°_{\overline{\eta}}, \mathscr{L}) = \mathrm{H}^i(Y_s, \Psi_h j_{\eta!}\mathscr{L})$ est unipotente. Il suffit de montrer que pour chaque $\alpha$, l'action de $\pi_1(\eta, \overline{\eta})$ sur $\Psi^\alpha_h(j_{\eta!}\mathscr{L})$ est unipotente, c'est-à-dire (**1.2.4**) qu'il existe une filtration de ce faisceau à gradué muni de l'action triviale. D'après **1.2.5**, le problème est local (pour la topologie étale) sur $Y_s$ (donc sur $Y$) et on peut alors supposer, par dévissage, que le faisceau localement unipotent $\mathscr{L}$ est *constant*, de valeur $\mathbb{Z}/n\mathbb{Z}$ pour un entier $n$ inversible sur $X$. Le triangle

$$\Psi_h(j_{\eta!}\mathbb{Z}/n\mathbb{Z}) \to \Psi_h(\mathbb{Z}/n\mathbb{Z}) \to \Psi_h(i_{\eta\star}\mathbb{Z}/n\mathbb{Z}) = i_{s\star}\Psi_{D/X}(\mathbb{Z}/n\mathbb{Z}) \to +1,$$

où $i$ est l'immersion fermée du diviseur $D$ dans $Y$ (et l'égalité résulte de [SGA 7 XIII, 1.3.6.1]), nous ramène à démontrer l'unipotence des faisceaux $\Psi^\alpha_h(\mathbb{Z}/n\mathbb{Z})$. Or, il est bien connu que l'action de l'inertie sur ces faisceaux est *triviale* (voir par exemple [Illusie 2004, §3.3] et [SGA 7 I, 3.3]).

**2.3.5**. *Unipotence : cas général.* On procède par réduction au cas précédent. Le théorème de changement de base propre nous permet d'une part de supposer le schéma de base $X$ strictement local intègre, de point fermé $x$ et de point *générique* $\eta \in X°$ — dont on note $\overline{\eta} \to \eta$ un point générique géométrique — et d'autre part d'exprimer la cohomologie $\mathrm{R}\Gamma(Y_{\overline{\eta}}, j_!\mathscr{L})$ de la fibre générique géométrique comme la cohomologie $\mathrm{R}\Gamma(Y_x, \Psi_{x,\overline{\eta}}(j_!\mathscr{L}))$ de la fibre spéciale $Y_x$ à coefficients dans un complexe naturellement muni d'une action de $\pi_1(\eta, \overline{\eta})$, calculant la cohomologie des fibres de Milnor :

$$\Psi_{x,\overline{\eta}} := (Y_x \hookrightarrow Y)^\star \mathrm{R}(Y_{\overline{\eta}} \to Y)_\star.$$

(Sur la théorie de P. Deligne des cycles évanescents *en dimension* $\geqslant 1$, voir par exemple [Laumon 1983].) D'après **1.2.3** et la normalité de $X°$, il suffit de montrer que l'action du groupe de Galois $\pi_1(\eta, \overline{\eta})$ sur chaque $\mathrm{H}^i_c(Y°_{\overline{\eta}}, \mathscr{L})$ est unipotente. La description précédente et le critère **1.2.5** nous ramènent à démontrer l'énoncé *local* suivant : pour tout entier $0 \leqslant \alpha \leqslant 2$ et tout point géométrique $\overline{y}$ de $Y$ localisé en $y \in (Y \smallsetminus Y°)_x$, l'action du groupe de Galois $\pi_1(\eta, \overline{\eta})$ sur le groupe de cohomologie $\mathrm{H}^\alpha(\overline{V}, j_!\mathscr{L}) = \Psi^\alpha_{x,\overline{\eta}}(j_!\mathscr{L})_{\overline{y}}$ de la *fibre de Milnor* $\overline{V} := Y_{(\overline{y})} \times_X \overline{\eta}$ est *unipotente*. (Rappelons que $j$ est l'immersion ouverte $Y° \hookrightarrow Y$.) Par hypothèse, la restriction du faisceau $\mathscr{L}$ à $V := Y_{(\overline{y})} \times_X \eta$ est unipotente ; par dévissage, on peut donc supposer le faisceau $\mathscr{L}$ constant de valeur $L$. Le triangle $j_!L \to L \to i_\star L \overset{+1}{\to}$ sur $Y_\eta$ nous ramène à montrer que l'action de $\pi_1(\eta, \overline{\eta})$ sur $\mathrm{H}^\alpha(\overline{V}, L)$ et $\mathrm{H}^{\alpha-1}(\overline{V} \times_Y D, L)$ est unipotente, où $D$ est le fermé de $Y$ constitué des diviseurs horizontaux. Par locale acyclicité des morphismes lisses seul est à considérer le premier groupe, dans le cas particulier où $y$ est un point singulier du morphisme $h$. Notons que le lieu singulier du morphisme $h$ est fini sur la base : la proposition [ÉGA II, 7.1.7] nous permet comme dans le lemme ci-dessus (et [Orgogozo 2006, §7.1]) de supposer que $X$ est un trait, par comparaison avec la théorie classique des cycles évanescents de Grothendieck. (On utilise ici un théorème de P. Deligne de commutation au changement de base, rappelé en [ibid., prop. 6.1].)

2.3.6. *Modération : cas général (esquisse).* Comme dans le paragraphe précédent, on se ramène à démontrer l'énoncé local suivant (pour $X$ strictement local, de point fermé $x$) : pour tout point géométrique $\overline{y}$ de $Y$ localisé en $y \in (Y \smallsetminus Y^{\circ})_x$, l'action du groupe de Galois $\pi_1(\eta, \overline{\eta})$ sur les groupes de cohomologie $\mathrm{H}^{\alpha}(\overline{V}, j_!\mathscr{L})$ de la *fibre de Milnor* $\overline{V} \coloneqq Y_{(\overline{y})} \times_X \overline{\eta}$ est d'ordre premier à l'exposant caractéristique $p$ de $x$. Nous traitons ici le cas où $Y \to X$ n'est pas lisse en $y$, l'autre cas étant semblable et plus élémentaire. D'après [F. Kato 2000, §2, p. 227][⑥] (voir aussi [Mochizuki 1995, §3]), le morphisme $h : Y \to X$ est sous-jacent à un morphisme *log-lisse* saturé $(Y, \mathscr{N}) \to (X, \mathscr{M})$ entre log-schémas fins et saturés (fs) tel que le lieu de trivialité de $(Y, \mathscr{N})$ (resp. de $(X, \mathscr{M})$) contienne la fibre $Y_\eta$ (resp. $\eta$). Notons $V$ la fibre $Y_{(\overline{y})} \times_X \eta$, normale par [Gabber et Ramero 2016, 10.7.17 (i)] ; la restriction à $V$ du faisceau $\mathscr{L}$ correspond à une représentation de $\pi_1(V)$, qui — par hypothèse de modération — se factorise à travers son plus grand quotient $\pi_1(V)_{\perp p}$ premier à $p$ : le groupe fondamental du schéma $V$ est un quotient du groupe de Galois du point générique $v$ de $V$ et l'action de ce dernier se factorise par hypothèse à travers un quotient premier à $p$. Les revêtements d'ordre premier à $p$ de $\overline{V}$ étant modérés, ils proviennent de la log-structure d'après le théorème [ibid., 11.4.39] d'acyclicité des morphismes log-lisses : il existe donc un morphisme Kummer-étale $(Y', \mathscr{N}') \to (Y, \mathscr{N})$ tel que la restriction de $\mathscr{L}$ à la fibre $V' = Y'_{(\overline{y})} \times_X \eta$ correspondante soit l'image inverse d'un faisceau sur $\eta$, de monodromie première à $p$. (On utilise ici le fait que le groupe fondamental de $V$ est extension de $\pi_1(\eta, \overline{\eta})$ par le groupe fondamental de $\overline{V}$.) Ceci nous ramène à démontrer la modération dans le cas particulier où $\mathscr{L}$ est un faisceau *constant*. Par commutation des cycles proches aux changements de base dans ce cas (rappelée en 2.3.5 *supra*), on peut supposer que $X$ est un *trait*, de corps des fractions $\kappa(\eta)$. La modération est connue dans ce cas (2.3.3).

**2.4. Image directe : immersion ouverte régulière.** Ce paragraphe est consacré à la démonstration de la proposition suivante, dont on trouvera une variante en [Gabber et Loeser 1996, lemme 4.3.2], et notamment d'un complément sur la taille des images directes.

**2.4.1. Proposition.** *Soient $j : X^{\circ} \hookrightarrow X$ une immersion ouverte telle que la paire $(X, X^{\circ})$ soit strictement régulière (2.1.1) et $\mathscr{L}$ un faisceau lisse de $\mathbb{Z}/n\mathbb{Z}$-modules sur $X^{\circ}$, où $n$ est un entier inversible sur $X$, tel que $j_!\mathscr{L}$ soit* modéré *(resp. localement unipotent le long de $\{X^{\circ}, X \smallsetminus X^{\circ}\}$). Soit $\mathfrak{D}$ la stratification de $X$ définie par la fonction semi-continue supérieurement associant à un point le nombre de branches de $D \coloneqq X \smallsetminus X^{\circ}$ y passant. Alors, les faisceaux $\mathrm{R}^c j_\star \mathscr{L}$ ($c \in \mathbb{N}$) sont constructibles et modérés (resp. localement unipotents) le long de la stratification $\mathfrak{D}$.*

*De plus :*

  (**i**) *si $X$ est lisse sur un schéma $S$ et $D = X \smallsetminus X^{\circ}$ est le complémentaire d'un diviseur à croisements normaux relativement à $S$, la formation de $\mathrm{R}^c j_\star \mathscr{L}$ commute aux changements de bases $S' \to S$ ;*

---

[⑥] Je remercie Arthur Ogus de m'avoir communiqué cette référence.

(ii) *si* $k : X \hookrightarrow \overline{X}$ *est une immersion ouverte et* $k_!j_!\mathscr{L}$ *modéré (resp. localement unipotent le long de la stratification* $\{X°, \overline{X} \smallsetminus X°\}$ *de* $\overline{X}$*), le complexe* $k_!\mathrm{R}j_\star\mathscr{L}$ *est à cohomologie* modérée *(resp. localement unipotente le long de la stratification définie par* $\mathfrak{D}$*, et le fermé* $F = \overline{X} \smallsetminus X$*).*

« Nombre de branches » : c'est le rang des fibres de $\mathrm{R}^1 j_\star \mathbb{F}_\ell$, pour $\ell$ inversible sur $X$, ou encore le rang de la log-structure naturellement associée à $(X, X°)$ — de sorte que $\mathfrak{D}$ est la *stratification canonique* ([STG VI, 1.5]) — : on choisit donc $\mathfrak{D}$ de sorte que le théorème soit automatiquement vrai pour les faisceaux constants.

Traitons tout d'abord l'assertion principale dans le cas constructible modéré ; la commutation au changement de base (i) est établie en 2.4.3 et le léger renforcement (ii) en 2.4.4. Le cas unipotent, plus facile, est traité en 2.4.5, directement sous la forme forte (ii).

**2.4.2**. *Cas modéré.* L'énoncé étant local pour la topologie étale, on peut supposer $X$ strictement local. Le diviseur $D$ est donc réunion transverse de diviseurs réguliers $D_i = \mathrm{div}(t_i)$, où les $t_1, \ldots, t_b$ sont des éléments de $\Gamma(X, \mathscr{O}_X)$. Comme dans [Deligne 1980, §1.7.8], posons $X_r := X[(t_i)^{1/r} : 1 \leq i \leq b]$ pour chaque entier $r \geq 1$. C'est un schéma régulier ([SGA 1 XIII, 4.1]) et $X_r° := X° \times_X X_r$ est le complémentaire d'un diviseur strictement à croisements normaux $D_r$. Notons $\pi_r$ le morphisme fini $X_r \to X$, $\pi_r° : X_r° \to X°$ le morphisme étale qui s'en déduit par restriction à l'ouvert $X°$ et $j_r : X_r° \hookrightarrow X_r$ l'immersion ouverte. L'hypothèse de modération assure que pour chaque $\mathscr{L}$ comme dans l'énoncé il existe un entier $r$, premier à la caractéristique résiduelle de $X$, tel que $\mathscr{L}_r := \pi_r^{°\star}\mathscr{L}$ soit un faisceau *constant* sur $X_r°$, de fibre que nous notons $L$ (droit).

$$\begin{array}{ccc}
\underline{L} = \mathscr{L}_r & X_r° \xrightarrow{\;j_r\;} X_r \\
& \pi_r° \downarrow \quad \square \quad \downarrow \pi_r \\
\mathscr{L} \hookrightarrow \mathscr{L}' & X° \xrightarrow{\;j\;} X
\end{array}$$

Le morphisme d'adjonction $\mathscr{L} \to \mathscr{L}' := \pi_{r\star}°\mathscr{L}_r$ étant injectif de conoyau lisse et modéré le long de $D$, il suffit — par récurrence sur le degré cohomologique — d'établir la proposition pour le faisceau $\mathscr{L}'$. (Comparer avec [Th. finitude, A, 1.3.3] et 2.3.3.) Par transitivité des images directes, on a $\mathrm{R}j_\star \mathscr{L}' = \pi_{r\star} \mathrm{R}j_{r\star}\mathscr{L}_r$. Puisque l'image directe par $\pi_r$ transforme un faisceau constructible modéré le long de $\mathfrak{D}_r$ en un faisceau constructible modéré le long de $\mathfrak{D}$, il suffit de traiter le cas $r = 1$, c'est-à-dire $\mathscr{L}$ constant. Cela résulte de la pureté cohomologique absolue, démontrée par O. Gabber (voir par exemple [STG XVI, 3.1.4]).

**2.4.3**. *(i) : commutation au changement de base.* La propreté cohomologique du morphisme *fini* $\pi_{r\star}$ nous permet, comme dans le paragraphe précédent, de nous ramener au cas particulier où le faisceau $\mathscr{L}$ est constant : on veut montrer que le faisceau $\mathrm{R}^c j_\star L$ commute aux changements de base $S' \to S$ lorsque $L$ est constant, fini, de torsion inversible sur $S$. Cela résulte de [Th. finitude, A, 1.3.3 (i)].

**2.4.4**. *(ii) : étude à l'infini, cas modéré.* On démontre maintenant l'énoncé 2.4.1 (ii) dans le cas modéré, en « passant à la limite sur $r$ ». (Cette présentation, plus élégante que que l'originale, nous a été proposée par O. Gabber.[7])

Fixons un entier $c \geqslant 0$. Le problème étant à nouveau local pour la topologie étale, on peut supposer $\overline{X}$ strictement local, d'exposant caractéristique $p \geqslant 1$. Soit $\overline{d}$ un point géométrique localisé en $d \in D \coloneqq X \smallsetminus X^\circ$, d'exposant caractéristique $p_d$ divisant $p$. (L'énoncé à démontrer est trivial si on se place en un point de $\overline{X} \smallsetminus X$ ou de $X^\circ$.) On veut montrer que l'action du groupe de Galois $\pi_1(d, \overline{d})$ sur $\mathrm{H}^c(X^\circ_{(\overline{d})}, \mathscr{L})$ — où $X^\circ_{(\overline{d})}$ est l'ouvert normal $X^\circ \times_X X_{(\overline{d})}$ de l'hensélisé strict en $\overline{d}$ — se factorise à travers un quotient premier à $p$. Par pureté, les dévissages du paragraphe précédent montrent que ce groupe de cohomologie est isomorphe à $\mathrm{H}^c(\pi_1^{\mathrm{mod}}(X^\circ_{(\overline{d})}, \overline{\eta}), L)$, où $L$ est la fibre de $\mathscr{L}$ en un point générique géométrique $\overline{\eta}$ de $X^\circ_{(\overline{d})}$ localisé en un point $\eta$ de $X^\circ_{(\overline{d})}$, et $\pi_1^{\mathrm{mod}}(X^\circ_{(\overline{d})}, \overline{\eta})$ est (non canoniquement) isomorphe à $\widehat{\mathbb{Z}}^{(p'_d)}(1)^b$, où $b$ est le nombre de branches de $D$ en $\overline{d}$. Le faisceau $\mathscr{L}$ étant modéré sur $\overline{X}$, l'action du quotient $\pi_1(X^\circ_{(d)}, \overline{\eta})$ de $\pi_1(\eta, \overline{\eta})$ sur $L$ se fait à travers un quotient d'ordre premier à $p$. La conclusion résulte alors du lemme ci-dessous, appliqué à la suite exacte déduite de $1 \to \pi_1(X^\circ_{(\overline{d})}) \to \pi_1(X^\circ_{(d)}) \to \pi_1(d, \overline{d}) \to 1$ par le quotient *caractéristique* $\pi_1(X^\circ(\overline{d})) \twoheadrightarrow \pi_1^{\mathrm{mod}}(X^\circ_{(\overline{d})}) = \pi_1(X^\circ_{(\overline{d})})/N$, où $G$ est le groupe de Galois, $A$ le groupe fondamental modéré $\pi_1^{\mathrm{mod}}(X^\circ_{(\overline{d})})$ et enfin $\pi$ le quotient $\pi_1(X^\circ_{(d)})/N$. Notons que les groupes d'inertie modérés en les points maximaux de $D$ engendrent $A$ et sont stables sous l'action naturelle de $G$ car les branches sont définies sur l'hensélisé (non strict) $X^\circ_{(d)}$.

**Lemme.** *Soient $p$ un nombre premier, $n \perp p$ un entier, $1 \to A \to \pi \to G \to 1$ une suite exacte de groupes profinis et $L$ un $\mathbb{Z}/n\mathbb{Z}$-module fini muni d'une action continue de $\pi$ se factorisant à travers un quotient premier à $p$. On suppose de plus le groupe $A$ abélien et produit fini, comme $G$-module, de groupes procycliques sans torsion d'ordre premier à $p$. Alors, pour tout entier $c \geqslant 0$, l'action de $G$ sur $\mathrm{H}^c(A, L)$ se factorise à travers $G_{\perp p}$.*

*Démonstration.* La décomposition $G$-équivariante de $A$ en produit, la suite spectrale de Hochschild-Serre et l'hypothèse $n \perp p$ — garantissant la stabilité par extension de la propriété à démontrer — nous ramènent immédiatement au cas particulier où $A$ est procyclique ; géométriquement, cela revient à supposer qu'il n'y a qu'une branche ($b = 1$ avec les notations précédentes). Pour un tel $A$, il est bien connu que la cohomologie de $L$ est concentrée en degrés $\{0, 1\}$ et $\mathrm{H}^0(A, L)$ (resp. $\mathrm{H}^1(A, L(1))$) sont respectivement les invariants et coïnvariants de $A \curvearrowright L$. (Voir par exemple [STG XVII, §1.2] pour une discussion de ces isomorphismes, indépendants du choix d'un progénérateur.) L'action de $\pi$ sur $L$ se factorisant par un quotient d'ordre premier à $p$, il en est de même de l'action de $G$ sur ces sous-quotients. □

**2.4.5**. *(ii) : cas unipotent.* Comme précédemment, on peut supposer $\overline{X}$ strictement local et, par hypothèse, $\mathscr{L}$ extension de faisceaux constants sur $X^\circ$ puis, par dévissage immédiat, constant et même isomorphe à $\mathbb{Z}/n\mathbb{Z}$ quitte à remplacer $n$ par un diviseur. Il s'agit donc de montrer la locale unipotence des faisceaux $k_!\mathrm{R}^c j_\star \mathbb{Z}/n\mathbb{Z}$, pour $n$ inversible sur $\overline{X}$. Cela résulte immédiatement de la pureté cohomologique absolue [STG XVI, 3.1.4].

---

[7] Ci-dessous, comme ailleurs, l'auteur est bien sûr seul responsable des erreurs.

Ceci achève la démonstration de la proposition 2.4.1. Nous souhaitons maintenant borner la taille des images directes par l'immersion ouverte $j$. (L'étude des fibres de $k_!\mathrm{R}j_\star$ s'y ramène immédiatement.)

**2.4.6. Proposition.** *Soit $c \geq 0$ un entier. Il existe une fonction $\varphi_c : \mathbb{N} \to \mathbb{N}$ telle que pour toute paire régulière et strictement locale $(X, X^\circ)$, tout entier $n \geq 1$ inversible sur $X$ et tout faisceau lisse $\mathscr{L}$ de $\mathbb{Z}/n\mathbb{Z}$-modules sur $X^\circ$ constructible et modéré le long de $X \smallsetminus X^\circ = \sum_{i=1}^b D_i$ le groupe de cohomologie $\mathrm{H}^c(X^\circ, \mathscr{L})$ soit $\varphi_c(b)$-extension ponctuelle de sous-quotients de $\mathscr{L}$, où $b$ désigne le nombre de branches de $D := X \smallsetminus X^\circ$.*

« extension ponctuelle » : cf. 1.7.

*Démonstration.* Reprenons les notations du paragraphe 2.4.2. Soient $r$ tel que $\pi_r^{\circ\star}\mathscr{L}$ soit constant, de fibre $L$, et $\boldsymbol{\mu}_{r,X}$ le groupe cyclique d'ordre $r$ des racines de l'unité d'ordre divisant $r$ dans $\Gamma(X, \mathscr{O}_X)$. Le revêtement $X_r^\circ \to X^\circ$ étant galoisien de groupe $\boldsymbol{\mu}_{r,X}^b$, le complexe $\mathrm{R}\Gamma(X^\circ, \mathscr{L})$ se réécrit $\mathrm{R}\Gamma(\boldsymbol{\mu}_{r,X}^b, \mathrm{R}\Gamma(X_r^\circ, L))$. Rappelons ([Serre 1968, VIII.§4]) que si $G = \langle \sigma \rangle$ est un groupe cyclique et $M$ un $G$-module, la cohomologie de $G$ à valeurs dans $M$ est représentée par le complexe périodique

$$0 \longrightarrow M \xrightarrow{\sigma - 1} M \xrightarrow{t} M \xrightarrow{\sigma - 1} M \xrightarrow{t} \cdots$$

où $t = \sum_{g \in G} g$, et le premier $M$ est placé en degré 0. Cette description, uniforme en $r$, nous ramène au cas particulier où $r = 1$, c'est-à-dire où $\mathscr{L}$ est constant, que l'on peut alors supposer isomorphe à $\mathbb{Z}/n\mathbb{Z}$. Dans ce cas, la conclusion résulte de la pureté : $\mathrm{H}^c(X^\circ, \mathbb{Z}/n\mathbb{Z}) = \bigwedge^c (\mathbb{Z}^b/n\mathbb{Z}^b)$.

Alternativement, on aurait pu passer à la limite sur $r$ et utiliser la description de $\mathrm{H}^c(A, L)$ faite en 2.4.4. Ceci permet d'ailleurs d'établir immédiatement l'annulation $\mathrm{H}^c(X^\circ, \mathscr{L}) = 0$ pour $c > b$. □

## 3. Image directe par un morphisme propre

Nous allons montrer que la propriété d'une famille d'être constructible et modérée (resp. localement unipotente) pour la topologie des altérations est stable par image directe propre. Cette stabilité est mise en défaut si l'on considère la topologie finie au lieu de la topologie des altérations (voir 3.2).

Ces résultats sont connus de longue date d'Ofer Gabber.

### 3.1. Énoncés.

**3.1.1. Théorème** (d'uniformité propre). *Soient $S$ un schéma nœthérien quasi-excellent, $f : X \to S$ un morphisme propre, $\alpha : X' \twoheadrightarrow X$ une altération et $\mathfrak{X}$ une stratification de $X'$. Il existe une altération $\beta : S' \twoheadrightarrow S$ et une stratification $\mathfrak{S}$ de $S'$ telles que pour tout entier $n \geq 1$, tout $i \geq 0$ et tout faisceau $\mathscr{F}$ de $\mathbb{Z}/n\mathbb{Z}$-modules sur $X[1/n]$ dont l'image inverse $\alpha[1/n]^\star \mathscr{F}$ est constructible et* *modérée (resp. constructible et localement unipotente) le long de $\mathfrak{X}[1/n]$ (au sens de* 1.1.3, 1.3.1, 1.2.1*) alors chaque $\beta[1/n]^\star \mathrm{R}^i f[1/n]_\star \mathscr{F}$ est constructible et modéré (resp. constructible et localement unipotent) le long de $\mathfrak{S}[1/n]$.*

$$\begin{array}{ccc} \mathscr{F} & & \mathfrak{X} \\ X & \xleftarrow{\alpha} & X' \\ f \downarrow & & \\ S & \xleftarrow{\exists \beta} & S' \\ & & \exists \mathfrak{S} \end{array}$$

*De plus, il existe un entier $N$, dépendant uniquement du triplet $(f, \alpha, \mathfrak{X})$, tel que chaque $\mathrm{R}^i f[1/n]_\star \mathscr{F}$ comme ci-dessus soit $N$-extension ponctuelle de sous-quotients de $\mathscr{F}$ relativement à $f$ (au sens de* 1.7*).*

**3.1.2.** Tout faisceau abélien constructible sur un schéma nœthérien étant modéré par un morphisme *fini* surjectif (cf. 1.6.7), on peut voir l'énoncé 3.1.1 comme une amélioration du théorème classique de constructibilité des images directes par un morphisme propre ([SGA 4 XIV, 1.1]), du moins lorsque l'on suppose la base nœthérienne quasi-excellente. (L'hypothèse d'excellence est d'ailleurs superflue ; voir 3.5.5 pour une esquisse d'argument, dont les détails sont laissés au lecteur.)

Il est également naturel d'exprimer le théorème précédent sous forme « dérivée ».

**3.1.3.** Soient $\alpha, \mathfrak{X}$ comme ci-dessus et $n$ inversible sur $X$. Notons $D^{b\,\mathrm{mod}}_{\alpha,\mathfrak{X}}(X, \mathbb{Z}/n\mathbb{Z})$ la sous-catégorie triangulée de $D^b_c(X, \mathbb{Z}/n\mathbb{Z})$ des complexes à cohomologie constructible le long de $\mathfrak{X}$ et modérée par $\alpha$. Il existe $\beta, \mathfrak{S}$ tels que $\mathrm{R}f_\star : D^b_c(X, \mathbb{Z}/n\mathbb{Z}) \to D^b_c(S, \mathbb{Z}/n\mathbb{Z})$ induise un foncteur entre les sous-catégories $D^{b\,\mathrm{mod}}_{\alpha,\mathfrak{X}}(X, \mathbb{Z}/n\mathbb{Z}) \to D^{b\,\mathrm{mod}}_{\beta,\mathfrak{S}}(S, \mathbb{Z}/n\mathbb{Z})$. (Comparer avec [Beilinson, Bernstein et Deligne 1982, §2.2.10-17].) Soit maintenant $\underline{n} = (n_i)_{i \in I}$ une famille d'entiers. Notons $D^{b\,\mathrm{mod}}_{\alpha,\mathfrak{X}}(X, \mathbb{Z}/\underline{n}\mathbb{Z})$ la catégorie triangulée produit $\prod_i D^{b\,\mathrm{mod}}_{\alpha,\mathfrak{X}}(X[1/n_i], \mathbb{Z}/n_i\mathbb{Z})$ et $D^{b\,\mathrm{mod}}_c(X, \mathbb{Z}/\underline{n}\mathbb{Z})_h$ la colimite sur $(\alpha, \mathfrak{X})$, des catégories triangulées $D^{b\,\mathrm{mod}}_{\alpha,\mathfrak{X}}(X, \mathbb{Z}/\underline{n}\mathbb{Z})$.

Si $f : X \to S$ est un morphisme propre, le théorème ci-dessus affirme que les foncteurs $\mathrm{R}f_\star : D^b_c(X[1/n_i], \mathbb{Z}/n_i\mathbb{Z}) \to D^b_c(S[1/n_i], \mathbb{Z}/n_i\mathbb{Z})$ induisent un foncteur

$$D^{b\,\mathrm{mod}}_c(X, \mathbb{Z}/\underline{n}\mathbb{Z})_h \to D^{b\,\mathrm{mod}}_c(S, \mathbb{Z}/\underline{n}\mathbb{Z})_h.$$

Mêmes énoncés pour les catégories $D^{b\,\mathrm{unip}}_{\alpha,\mathfrak{X}}(X, \mathbb{Z}/n\mathbb{Z})$ et $D^{b\,\mathrm{unip}}_c(X, \mathbb{Z}/\underline{n}\mathbb{Z})_h$, définies de manière évidente.

Ces énoncés sont non triviaux même lorsque tous les $n_i$ sont égaux : si l'on ne considère par exemple que des $\mathbb{F}_\ell$-faisceaux dans le théorème 3.1.1 (pour un $\ell$ fixé), la conclusion ne résulte pas directement, pour $I$ de cardinal $> 1$, du théorème de constructibilité usuel et fixer un $\ell$ ne semble d'ailleurs pas simplifier la démonstration. Notons cependant que si les entiers $n_i$ sont des puissances d'un même nombre premier $\ell$, la généralité apportée par l'ensemble d'indices $I$ est illusoire : un $\mathbb{Z}/\ell^r\mathbb{Z}$-faisceau constructible $\mathscr{F}$ est muni d'une filtration $\mathscr{F}_\alpha := \ell^\alpha \mathscr{F}$ de gradués des $\mathbb{F}_\ell$-faisceaux constructibles, lisses (resp. modérés, localement unipotents) si et seulement si il en est ainsi de $\mathscr{F}$. (Si $\mathscr{F}$ est même plat, les gradués non nuls sont isomorphes à $\mathscr{F}/\ell\mathscr{F}$.)

**3.1.4**. **Remarque.** Soient $X$ un schéma propre et lisse sur un corps algébriquement clos $k$ de caractéristique $p \geqslant 0$, $D$ un diviseur à croisements normaux et $U$ l'ouvert complémentaire. Il résulte du théorème précédent qu'il existe une constante $C_{X,D}$ telle que pour tout $\ell \neq p$ et tout $\mathbb{F}_\ell$-faisceau lisse $\mathscr{L}$ sur $U$, modéré le long de $D$, on ait $h^i_c(U, \mathscr{L}) \leqslant C_{X,D} \cdot \mathrm{rang}(\mathscr{L})$. Comme l'a expliqué Claude Sabbah à l'auteur, on peut en donner la démonstration suivante lorsque $k$ est le corps $\mathbb{C}$ de nombres complexes. Soit $T$ un voisinage tubulaire *ouvert* de $D(\mathbb{C})$ tel que l'inclusion du *compact* $K = X(\mathbb{C}) - T$ dans $U(\mathbb{C})$ soit une équivalence d'homotopie ; le groupe de cohomologie $\mathrm{H}^i_{\mathrm{Betti},c}(U(\mathbb{C}), \mathscr{L})$ est donc isomorphe à $\mathrm{H}^i_{\mathrm{Betti},c}(K, \mathscr{L})$. Par compacité, il existe un recouvrement *fini* $V_i$ de $K$ par des ouverts simplement connexes. Il en résulte que la cohomologie de $U(\mathbb{C})$ à valeurs dans $\mathscr{L}$ est l'aboutissement d'une suite spectrale dont les objets de la première « page » ne dépendent que des $V_i$ et du rang de $\mathscr{L}$. On conclut par le théorème de comparaison Betti-étale d'Artin-Grothendieck. Il serait sans doute intéressant d'adapter cet argument au cas d'un corps quelconque à l'aide d'un voisinage tubulaire adéquat.

**3.2**. **Un exemple.** L'objet de cette section est de montrer que l'image directe d'une famille de faisceaux *constants* constructibles (donc, *a fortiori*, constructible et modérée pour la topologie grossière) n'est pas nécessairement modérée pour la topologie finie. Il est pour cela nécessaire de considérer une base de dimension $> 1$. On commence par rappeler une construction exprimant un groupe de monodromie locale, en dimension 2, comme un groupe de monodromie *globale* en dimension 1. (Ce résultat est bien connu ; cf. par exemple [Kato et Saito 2013, 2.6.3.1].) On rappelle ensuite un résultat de Igusa garantissant l'existence de familles à un paramètre à grosse monodromie.

**3.2.1**. Soient $k$ un corps algébriquement clos, $U_0$ un ouvert de $X_0 = \mathbf{P}^1_k$ et $X = \mathbf{A}^2_k$. On considère l'application rationnelle $X \dashrightarrow X_0$, $(u, v) \mapsto t = u/v$ et l'on note $U \subseteq X$ l'ouvert complémentaire des droites $u = \lambda v$ pour $\lambda \in X_0 - U_0$. Notons $O \in X$ l'origine, $X_{(O)}$ l'hensélisé (strict) et $U_{(O)}$ le produit fibré $X_{(O)} \times_X U$. Le morphisme $U \to U_0$ induit un morphisme entre les groupes fondamentaux, de même que $U_{(O)} \to U$. Par composition, on en déduit un morphisme $\pi_1(U_{(O)}) \to \pi_1(U_0)$ ; nous allons voir que ce dernier morphisme est *surjectif*. (La surjectivité de $\pi_1(U) \to \pi_1(U_0)$ est quant à elle triviale, $U \to U_0$ ayant une section.)

Soit $V_0$ un revêtement fini étale connexe de $U_0$ ; on veut montrer que le produit fibré $V_0 \times_{U_0} U_{(O)}$ est également connexe. Notons $Y_0$ la normalisation de $X_0$ dans $V_0$ et $Y$ celle de $X$ dans l'image inverse $V := V_0 \times_{U_0} U$. Montrons que la fibre de $Y \to X$ au-dessus de l'origine $O$ est connexe, c'est-à-dire est un point. Soit $X'$ l'éclatement de $X$ en l'origine ; c'est un fibré en $\mathbf{A}^1$ sur $X_0$. Il en est donc de même du produit fibré $Y' := X' \times_{X_0} Y_0$ sur $Y_0$. De la normalité de $Y'$, donc de $\mathcal{O}(Y')$, et de l'égalité $Y'_U = V$, on tire que le morphisme $Y \to X$ est la factorisation de Stein du morphisme $Y' \to X$. La fibre de $Y' \to X$ au-dessus de l'origine étant $Y_0$, donc connexe, il en est de même de celle de $Y \to X$[8].

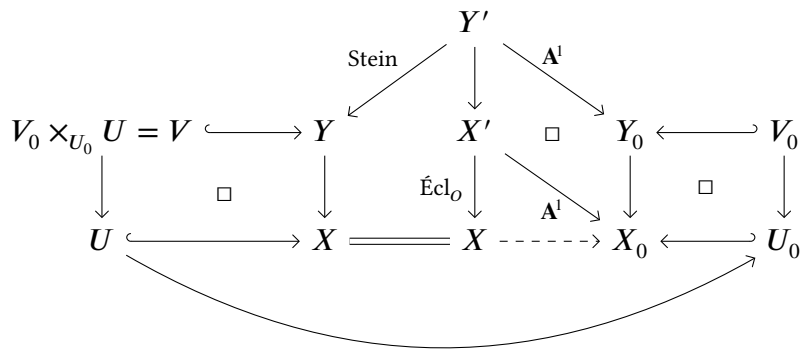

**3.2.2.** Soient maintenant $p$ un nombre premier impair, $k = \mathbb{F}_p$, $g \geqslant 1$ et $f \in k[x]$ un polynôme de degré $2g$ sans facteur carré. Notons $U_0$ l'ouvert $\mathrm{Spec}(k[t][f(t)^{-1}])$ de $\mathbf{A}_k^1$ et $\pi_{U_0} : \mathscr{C} \to U_0$ la courbe projective lisse d'équation affine $y^2 = f(x)(t-x)$. Pour chaque nombre premier impair $\ell \neq p$, notons $\mathscr{F}_\ell$ le $\mathbb{F}_\ell$-faisceau lisse $\mathrm{R}^1\pi_{U_0\star}\mathbb{F}_\ell$ sur $U_0$. D'après [Igusa 1959] lorsque $g = 1$ et [Yu 1995] en général, le groupe de monodromie géométrique de $\mathscr{F}_\ell$ est tout le groupe symplectique $\mathrm{Sp}(2g, \mathbb{F}_\ell)$. Rappelons rapidement un argument possible : le prolongement par zéro de $\mathscr{F}_\ell$ à $\mathbf{A}_k^1$ étant convolution additive de deux faisceaux de Kummer (quadratiques) explicites (voir [Katz et Sarnak 1999, 10.1.17]), les résultats généraux de [Katz 1988, chap. 8] entraînent que le faisceau $\mathscr{F}_\ell$ est « Fourier-irréductible » au sens d'*op. cit.* La conclusion résulte alors de la théorie de Picard-Lefschetz — la monodromie est engendrée par des transvections — et d'un théorème de théorie des groupes ([Zalesskiĭ et Serežkin 1976]). (Voir aussi [Katz et Sarnak 1999, p. 293–300] ainsi que [Deligne 1974a, 5.11] pour des variantes $\ell$-adiques.)

**3.2.3.** Soit $U$ l'ouvert de $X = \mathbf{A}_k^2$ défini à partir de $U_0$ comme en 3.2.1. On étend la courbe $\pi_U := \pi_{U_0} \times_{U_0} U$ en une courbe propre $\pi_X$ au-dessus de $X$ par normalisation, en voyant $\mathscr{C} \times_{U_0} U$ comme un revêtement double de $\mathbf{P}_U^1$. Pour chaque $\ell \neq p$ impair, la monodromie locale en l'origine $O \in X$ du faisceau lisse $\mathrm{R}^1\pi_{U\star}\mathbb{F}_\ell$ sur $U$ est, d'après les résultats des deux paragraphes précédents, le groupe $\mathrm{Sp}(2g, \mathbb{F}_\ell)$. Rappelons[9] que son cardinal est le produit $\ell^{g^2} \prod_1^g (\ell^{2i} - 1)$. Il en résulte que lorsqu'une puissance $p^r$ de $p$ divise $\ell - 1$ — ce qui se produit pour une infinité de nombres premiers $\ell$ — les $p$-Sylow de la monodromie locale en l'origine sont de cardinal au moins $p^r$.

Ceci a pour conséquence qu'il n'existe pas de $X$-morphisme fini surjectif modérant simultanément tous les faisceaux $\mathrm{R}^1\pi_{X\star}\mathbb{F}_\ell$, pour $\ell \neq p$ impair. Par contre, l'éclatement $X' \to X$ convient car chaque $\mathrm{R}^1\pi_{X_0\star}\mathbb{F}_\ell$ est modéré sur $X_0$, de même que son tiré en arrière par $X' \to X_0$.

---

[8]L'auteur remercie Takeshi Saito de son aide dans la rédaction de cet argument.

[9]Voir par exemple [E. Artin 1957, chap. III], ou [Serre 2007, 6.5.1] et [Steinberg 1968, §9, th. 25] pour des résultats plus généraux.

**3.3. Premières réductions.** Nous commençons la démonstration du théorème 3.1.1 — dont nous reprenons les notations — par des réductions diverses. Le complément sur la majoration des fibres est traité en 3.6.

**3.3.1.** *Réduction au cas où $X' = X$.* Nous vérifions ici que le théorème d'uniformité propre découle du cas particulier où l'altération $\alpha : X' \to X$ de la source est l'identité. Soit $X_\bullet$ le cosquelette de $\alpha$ — c'est un hyperrecouvrement propre de $X$ ([SGA 4 Vbis], [Deligne 1974b, 5.3.8]) — et notons $\mathfrak{X}_j$ ($j \geqslant 0$) les stratifications déduites de $\mathfrak{X}$ sur $X' = X_0$ par image inverse sur chaque $X_j$. Pour tout faisceau $\mathscr{F}$ sur $X[1/n]$ comme dans l'énoncé, on a $\mathrm{R}f[1/n]_\star \mathscr{F} = \mathrm{R}f_\bullet[1/n]_\star \mathscr{F}_\bullet$, où $f_\bullet$ est le morphisme $X_\bullet \to S$ et $\mathscr{F}_\bullet$ le tiré en arrière de $\mathscr{F}$ sur (le topos total du schéma simplicial) $X_\bullet[1/n]$.

Si $\mathscr{F}$ est constructible et localement unipotent le long de $\mathfrak{X}[1/n]$, il en est de même de chaque $\mathscr{F}_j$ le long des $\mathfrak{X}_j[1/n]$ (1.4). De même, si $\mathscr{F}$ est modéré par $\alpha[1/n]$, chaque faisceau $\mathscr{F}_j$ sur $X_j[1/n]$ est modéré.

Soit $d$ un majorant de la dimension des fibres de $f$. Supposons que pour chaque $j \leqslant 2d$, le théorème soit démontré pour le morphisme propre $f_j : X_j \to S$, l'altération identité et la stratification $\mathfrak{X}_j$ de $X_j$. Notons $\beta_j$, $\mathfrak{S}_j$ et $N_j$ respectivement les altérations de $S$, stratifications et entiers associés.

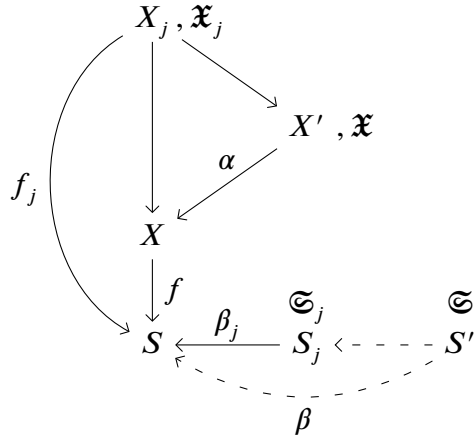

Par stabilité par image inverse des propriétés considérées on peut supposer les morphismes $\beta_j$ tous égaux à un même morphisme $\beta : S' \to S$. (L'ensemble partiellement ordonné des altérations de $S$ est filtrant [à droite].) Quitte à raffiner les stratifications $\mathfrak{S}_j$, on peut également supposer qu'elles sont toutes égales à une même stratification $\mathfrak{S}$ de $S'$. (Même argument.) D'après la suite spectrale de descente ([ibid., 5.2.7.1], [SGA 4 Vbis, 3.3.3])

$$E_1^{ji} = \mathrm{R}^i f_j[1/n]_\star \mathscr{F}_j \Rightarrow \mathrm{R}^{i+j} f[1/n]_\star \mathscr{F}$$

chaque $\mathrm{R}^k f[1/n]_\star \mathscr{F}$ est extension itérée de sous-quotients des $\mathrm{R}^i f_j[1/n]_\star \mathscr{F}_j$ pour $i+j = k$ et, par hypothèse, nul pour $k > 2d$. La constructibilité et locale unipotence (resp. modération) le long de $\mathfrak{S}$ des images directes par $f$ d'un faisceau comme dans l'énoncé sont alors conséquence de l'hypothèse sur $(\beta, \mathfrak{S})$ relativement aux $f_j$ et de la stabilité par extension de ces propriétés.

**3.3.2**. *Réduction au cas où la base $S$ est intègre.* Les réductions au cas connexe et réduit sont immédiates. On peut également supposer $S$ intègre car le coproduit de ses composantes irréductibles en est une altération.

**3.3.3**. *Réduction au cas où les faisceaux sont le prolongement par zéro de faisceaux lisses sur une même strate dense.* Soient $n \geq 1$ un entier et $\mathscr{F}$ un faisceau de $\mathbb{Z}/n\mathbb{Z}$-modules sur $X^\dagger := X[1/n]$ (notation de **1.5.2**). Si $\mathscr{F}$ est constructible le long de la stratification $\mathfrak{X}^\dagger$ — que l'on peut supposer bonne, quitte à la raffiner —, il est d'après **1.1.4** extension successive de faisceaux $i_\star^\dagger j_!^\dagger \mathscr{L}$, où :

- $j$ est l'immersion ouverte d'une strate $W$ de $\mathfrak{X}$ dans son adhérence $Y := \overline{W}$ ;
- $i$ est l'immersion fermée de $Y$ dans $X$ ;
- $\mathscr{L}$ est un faisceau abélien lisse de $\mathbb{Z}/n\mathbb{Z}$-modules sur $W^\dagger$.

$$\begin{array}{ccccc}
\mathscr{L} & & & & \mathscr{F} \\
W & \xhookrightarrow{j} & Y = \overline{W} & \xhookrightarrow{i} & X \\
& & \searrow_{g := f \circ i} & & \downarrow f \\
& & & & S
\end{array}$$

Les faisceaux $i_\star^\dagger j_!^\dagger \mathscr{L}$ sont constructibles (resp. modérés, localement unipotents) le long de $\mathfrak{X}^\dagger$.

Réciproquement, la stabilité par extension de ces mêmes propriétés et la suite exacte longue de cohomologie associée à une extension de faisceaux nous ramènent à démontrer le théorème pour les faisceaux $\mathscr{F} = i_\star^\dagger j_!^\dagger \mathscr{L}$ comme ci-dessus. Par finitude du nombre de strates et la propriété des altérations d'être filtrantes à droite (cf. **3.3.1**), on peut supposer les immersions $i$ et $j$ fixées. Posons $\mathfrak{Y} := i^\star \mathfrak{X}$, la stratification déduite de $\mathfrak{X}$ par restriction au fermé $Y$, et $g := f \circ i : Y \to S$ le morphisme composé. Le faisceau $j_!^\dagger \mathscr{L} = i^{\dagger\star} i_\star^\dagger j_!^\dagger \mathscr{L}$ étant modéré (resp. localement unipotent) le long de $\mathfrak{Y}^\dagger$ si $i_\star^\dagger j_!^\dagger \mathscr{L}$ l'est le long de $\mathfrak{X}^\dagger$, il résulte de la formule $\mathrm{R}g_\star^\dagger(j_!^\dagger \mathscr{L}) = \mathrm{R}f_\star^\dagger(i_\star^\dagger j_!^\dagger \mathscr{L})$, que l'on peut faire l'hypothèse supplémentaire que $Y = X$, c'est-à-dire que $\mathscr{F}$ est de la forme $j_!^\dagger \mathscr{L}$.

## 3.4. Énoncé de la récurrence nœthérienne.

**3.4.1**. Pour chaque triplet $(S, \rho, d)$, où $S$ est un schéma nœthérien quasi-excellent intègre ou vide et $\rho \geq -1$, $d \geq -1$ sont des entiers, on considère les deux énoncés suivants :

> $\mathsf{Mod}(S, \rho, d)$ : pour toute altération intègre $T$ de $S$, tout $T$-schéma propre $f : Y \to T$ de fibre générique de dimension au plus $d$ et toute stratification $\mathfrak{Y}$ de $Y$, il existe une altération $\beta : T' \to T$, une stratification $\mathfrak{T}$ de $T$ et un entier $N$ tels que pour tout entier $n \geq 1$ et tout faisceau $\mathscr{F}$ de $\mathbb{Z}/n\mathbb{Z}$-modules sur $Y[1/n]$ constructible le long de $\mathfrak{Y}[1/n]$ et modéré, chaque image directe $\mathrm{R}^k f[1/n]_\star \mathscr{F}$ est, pour $k \leq \rho$, constructible le long de $\mathfrak{T}[1/n]$, modérée par $\beta[1/n]$ et $N$-extension ponctuelle de $\mathscr{F}$ relativement à $f[1/n]$ ;

et

> Unip$(S, \rho, d)$ : pour toute altération intègre $T$ de $S$, tout $T$-schéma propre $f : Y \to T$ de fibre générique de dimension au plus $d$ et toute stratification $\mathfrak{Y}$ de $Y$, il existe une altération $\beta : T' \to T$ et une stratification $\mathfrak{T}$ de $T'$ telles que pour tout entier $n \geqslant 1$ et tout faisceau $\mathscr{F}$ de $\mathbb{Z}/n\mathbb{Z}$-modules sur $Y[1/n]$ constructible et localement unipotent le long de $\mathfrak{Y}[1/n]$, chaque image directe $\mathrm{R}^k f[1/n]_\star \mathscr{F}$ est, pour $k \leqslant \rho$, constructible et localement unipotente le long de $\mathfrak{T}[1/n]$.

$$\begin{array}{c} Y\,\mathfrak{Y} \\ f\downarrow \\ S \xleftarrow{\text{altération}} T \xleftarrow{\exists \beta} T' \end{array}$$

Commençons par quelques remarques sur ces énoncés.

**3.4.2.** Le premier énoncé (Mod) affirme l'existence d'une stratification de $T$ et non de $T'$ : cela revient au même d'après le résultat de descente de stratification 1.5.4. Dans Unip, le complément sur les extensions ponctuelles n'est pas inclus ; d'après l'observation 1.3.3, cela résulte du cas modéré.

Si $T$ est une altération intègre d'un schéma $S$ comme ci-dessus, les propriétés Mod$(T, \rho, d)$ et Mod$(S, \rho, d)$ sont *équivalentes* ; de même pour la variante (localement) unipotente : on démontre une propriété de l'espace de Zariski-Riemann.

Enfin, notons que ces énoncés sont trivialement vrais pour $S = \emptyset$ ou $\rho = -1$ ; ce n'est pas le cas pour $d = -1$ (morphisme non dominant).

**3.4.3.** Fixons $(S, \rho, d)$ et précisons la différence entre l'énonce Mod$(S, \rho, d)$ ci-dessus et l'énoncé analogue où l'on considère plus généralement des faisceaux $\mathscr{F}$ de $\mathbb{Z}/n\mathbb{Z}$-modules constructibles modérés le long d'une stratification non pas de $Y$ mais d'une (même) altération $\alpha : Y' \to Y$. La formule de descente cohomologique $\mathrm{R}f_\star^\dagger \mathscr{F} = \mathrm{R}f_{\bullet\star}^\dagger \mathscr{F}_\bullet$, où la source de $f_\bullet$ est le 0-cosquelette de $\alpha$, ramène (comme en 3.3.1) cet énoncé *a priori* plus général à la conjonction de Mod$(S, \rho, d)$ et de Mod$(S, \rho - 1, \ast)$ c'est-à-dire des Mod$(S, \rho - 1, d')$ pour $d' \geqslant -1$.

Même résultat pour Unip, *mutatis mutandis*.

**3.4.4.** Nous allons démontrer le théorème 3.1.1 par récurrence nœthérienne en considérant la relation bien fondée (=nœthérienne) $\prec$ sur les triplets $(S, \rho, d)$ lexicographiquement déduite des relations d'ordre usuelles sur les entiers $\rho$ et $d$ et de la relation binaire : $S' \prec S$ si et seulement si $S'$ est (isomorphe à) un fermé strict d'une altération intègre de $S$.

Notons tout de suite qu'il résulte du « lemme de Gabber » (1.6.2) que la propriété Mod ou Unip d'un triplet $(S, \rho, -1)$ est conséquence de cette même propriété pour les triplets $(\ast \prec S, \rho, \ast)$.

**Lemme.** *La relation $S' \prec S$ précédente est bien fondée : il n'existe pas de suite infinie $(S_n)$ de schémas nœthériens intègres telle que $S_{n+1}$ soit un fermé strict d'une altération $S'_n$ de $S_n$.*

$$
\begin{array}{c}
S_{n+1} \hookrightarrow S'_n \\
\downarrow \quad\quad \downarrow \\
S_n \\
\vdots \\
S_1 \hookrightarrow S'_1 \\
\downarrow \\
\overline{\{\tau\}} = F_N \subseteq \cdots \subseteq F_{n+1} \subseteq F_n \quad \subseteq \quad S_0 := S
\end{array}
$$

*Démonstration.* On raisonne par l'absurde, en supposant qu'il existe une suite infinie $\cdots \to S_{n+1} \hookrightarrow S'_n \to S_n \hookrightarrow \cdots$ ($n \geq 0$), où chaque $S'_n \to S_n$ est une altération et chaque $S_{n+1} \hookrightarrow S'_n$ est une immersion fermée stricte. Les images $F_n$ des $S_n$ dans $S := S_0$ forment une suite décroissante, donc stationnaire, de fermés de $S$. Il existe donc un entier $N$ tel que les points génériques des schémas $S_n$, pour $n \geq N$, soient tous d'image un même point $\tau \in S$. Pour chaque tel $n$, la restriction de $S'_n \to S_n$ (resp. $S_{n+1} \hookrightarrow S'_n$) aux fibres au-dessus de $\tau$ est encore une altération (resp. une immersion fermée *stricte*). Ceci nous ramène au cas particulier où les $S_n$ sont des schémas algébriques sur un corps, auquel cas la conclusion résulte des inégalités $\dim(S_{n+1}) < \dim(S_n) < +\infty$ pour chaque $n \geq 0$. □

**3.4.5.** Il résulte des dévissages précédents (§3.3) et de la nullité des $R^i f_\star$ pour $i$ grand que la conjonction des énoncés $\mathsf{Mod}(S, \rho, d)$ et $\mathsf{Unip}(S, \rho, d)$ pour $S$ nœthérien intègre quasi-excellent, $\rho \geq -1$ et $d \geq -1$ entraîne le théorème d'uniformité propre 3.1.1. La fin de cette section 3 est donc consacrée à la démonstration par récurrence de ces énoncés. À cet effet, on fixe un schéma nœthérien $S$, des entiers $\rho \geq -1, d \geq -1$ — que l'on peut supposer $\geq 0$ par récurrence et le lemme de Gabber — et l'on considère une altération intègre $T$ de $S$, un morphisme propre (dominant) $f : Y \to T$ de fibre générique de dimension au plus $d$ et une stratification $\mathfrak{Y}$ de $Y$. D'après 3.3.3 (démonstration), on peut supposer que chaque faisceau test $\mathscr{F}$ des énoncés $\mathsf{Mod}(S, \rho, d)$ et $\mathsf{Unip}(S, \rho, d)$ est le prolongement par zéro d'un faisceau lisse de $\mathbb{Z}/n\mathbb{Z}$-modules sur la strate ouverte dense $Y°[1/n]$ de $\mathfrak{Y}[1/n]$, pour un entier $n \geq 1$ (dépendant du faisceau).

### 3.5. Réduction au cas des courbes.

**3.5.1.** *Morphismes plurinodaux ([DE JONG 1997, 5.8]) : définitions et ubiquité.* Soit $f : Y \to T$ un morphisme entre schémas intègres et soit $Y°$ (resp. $T°$) un ouvert dense de $Y$ (resp. $T$). On dit que $f$ est **plurinodal adapté à** $(Y°, T°)$ s'il existe :

— une factorisation
$$f = (f_d : Y_d \to Y_{d-1}) \circ \cdots \circ (f_1 : Y_1 \to Y_0),$$
où $Y_0 = T$ et $Y_d = Y$, $d \in \mathbb{N}$ ;
— des ouverts denses $Y°_i \subseteq Y_i$, où $Y°_0 = T°$ et $Y°_d = Y°$,

tels que chaque $f_i$ soit une courbe nodale adaptée à $(Y_i^\circ, Y_{i-1}^\circ)$ au sens de 2.1.2.

Un morphisme $f : Y \to T$ est **plurinodal** (resp. **plurinodal adapté** à un ouvert dense $Y^\circ$ de $Y$) s'il existe des ouverts denses $Y^\circ$ et $T^\circ$ (resp. un ouvert dense $T^\circ$) de sorte que $f$ soit adapté à $(Y^\circ, T^\circ)$.

Un morphisme non nécessairement dominant $f : Y \to T$ est dit **presque plurinodal** (resp. **presque plurinodal adapté à un ouvert dense $Y^\circ$ de $Y$**) s'il est *somme* de morphismes plurinodaux entre schémas intègres et d'un morphisme non dominant (resp. et si pour chaque composante connexe $Z$ de $Y$ dominant $T$, le morphisme composé $Z \to T$ est plurinodal adapté à $Z^\circ := Y^\circ \cap Z$).

L'ubiquité de ces morphismes est assurée par le théorème suivant ([DE JONG 1997, corollaire 5.10]).

**Théorème** (de Jong). *Soient $f : Y \to T$ un morphisme propre, où $T$ est un schéma intègre nœthérien quasi-excellent, et $Y^\circ$ un ouvert dense de $Y$. Il existe une altération $T' \to T$ et un hyperrecouvrement propre $Y'_\bullet \to Y \times_T T'$ tel que $Y'_0 \to T'$ soit presque plurinodal adapté à un ouvert $Y'^{\circ\circ}_0$ contenu dans $Y'^{\circ}_0 := Y^\circ \times_Y Y'_0$.*

(L'énoncé original est légèrement différent ; voir [ORGOGOZO 2006, §4.3] si nécessaire pour une discussion.)

Comme expliqué en [DE JONG 1997, 5.16], l'hypothèse d'excellence n'est pas nécessaire : on peut s'en affranchir par passage à la limite.

**3.5.2.** Reprenons les notations et hypothèses de la fin du paragraphe 3.4.5 et montrons que l'on peut supposer le morphisme $f : Y \to T$ plurinodal adapté à la strate dense. On peut supposer que l'altération $T' \to T$ du théorème précédent est l'identité car l'énoncé que l'on souhaite démontrer est local en bas pour la topologie des altérations ; il est même local en haut pour cette topologie, de façon compatible à la récurrence, cf. 3.4.3. La descente cohomologique (cf. 3.3.1) et l'hypothèse de récurrence pour les triplets $(S, \rho - 1, *)$ nous permettent alors de supposer de plus le morphisme $f$ *presque* plurinodal : $f = f_+ \sqcup f_-$, où $f_+$ est plurinodal adapté à un ouvert $Y^{\circ\circ}$ de la strate dense $Y^\circ \subseteq Y$ et $f_-$ est un morphisme propre non surjectif, d'image $F \subseteq T$. Il résulte de l'hypothèse de récurrence pour les triplets $(F, \rho, *)$ — avec $F \prec S$ — et du lemme de Gabber 1.6.2 que l'on peut supposer $f = f_+$.

Notons respectivement $j^\circ$ et $j^{\circ\circ}$ les immersions ouvertes de $Y^\circ$ et $Y^{\circ\circ}$ dans $Y$, et fixons un entier $n \geq 1$. Pour tout faisceau test $\mathscr{F} = j^{\circ\dagger}_! \mathscr{L}^\circ$, où $\mathscr{L}^\circ$ est un faisceau lisse de $\mathbb{Z}/n\mathbb{Z}$-modules sur $Y^{\circ\dagger}$, on a une suite exacte

$$0 \to j^{\circ\circ\dagger}_! \mathscr{L}^{\circ\circ} \to j^{\circ\dagger}_! \mathscr{L}^\circ \to \overline{\mathscr{L}} \to 0,$$

où $\mathscr{L}^{\circ\circ} := \mathscr{L}^\circ_{|Y^{\circ\circ\dagger}}$, et le conoyau $\overline{\mathscr{L}}$ est un faisceau constructible à support dans le fermé strict $Z^\dagger := (Y - Y^{\circ\circ})^\dagger$ de $Y$. Pour $\mathscr{L}^\circ$ variable — mais astreint aux conditions déduites de celles sur $\mathscr{F}$ comme ci-dessus —, les $\overline{\mathscr{L}}$ forment une famille (uniformément) constructible (quitte à raffiner $\mathfrak{Y}$) ; modérée (resp. localement unipotente) si la famille des faisceaux $j^{\circ\dagger}_! \mathscr{L}^\circ$ l'est : cela résulte de la stabilité par quotient (cf. 1.4). La suite exacte longue de cohomologie associée à la suite exacte précédente et la stabilité par extension, maintes fois utilisée, nous ramènent à considérer des faisceaux test de la forme $j^{\circ\circ\dagger}_! \mathscr{L}^{\circ\circ}$ ou bien $\overline{\mathscr{L}}$. Dans ce dernier cas, on peut utiliser l'hypothèse de récurrence pour le triplet $(S, \rho, d-1)$ car la dimension de la fibre générique de $Z$ sur $T$ est strictement inférieure à celle de $Y$. (On rappelle que l'on a supposé $d \geq 0$, c'est-à-dire $f$ dominant.)

Comme annoncé, on peut donc supposer — quitte à raffiner $\mathfrak{Y}$ à nouveau — que l'on a l'égalité $Y^{\circ\circ} = Y^\circ$, c'est-à-dire que le morphisme $f$ est plurinodal adapté à une strate dense $Y^\circ$ de $\mathfrak{Y}$. On peut également supposer $Y$ intègre.

**3.5.3.** Sous les hypothèses en vigueur, on peut factoriser le morphisme $f : Y \to T$ en $f = g \circ h$, où $h : (Y, Y^\circ) \to (X, X^\circ)$ est une courbe nodale adaptée et $g : (X, X^\circ) \to (T, T^\circ)$ est un morphisme plurinodal adapté.

$$\begin{array}{c} Y \\ \downarrow h \\ f \quad X \\ \downarrow g \\ T \end{array}$$

D'après le théorème de résolution par altération des singularités des schémas quasi-excellents d'O. Gabber (voir [STG VII, 1.1] et [STG II, 3.2.1]), il existe des diagrammes commutatifs indicés par $i$

$$\begin{array}{ccc} T & \longleftarrow & T' \quad \text{régulier} \\ \text{altération} \uparrow & & \uparrow \\ \widetilde{T} & \longleftarrow & T'_i \end{array}$$

où $\widetilde{T}$ est une même altération de $T$, les $T'_i$ forment un recouvrement fini de Zariski de $\widetilde{T}$, et la paire $(T', T'^\circ := T' \times_T T^\circ)$ est *régulière* (2.1). D'après 2.3.1, appliqué au morphisme $h' := h \times_T T' : Y' \to X'$ et aux ouverts $Y'^\circ := Y^\circ \times_T T'$, $X'^\circ := X^\circ \times_T T'$, pour chaque $k \in \{0, 1, 2\}$ la famille des $R^k h'^\dagger_\star(\mathscr{F}')$ est constructible et modérée (resp. localement unipotente) sur $X'$, lorsque $\mathscr{F}$ parcourt les faisceaux test du type considéré : prolongement par zéro d'un faisceau lisse, modéré (resp. localement unipotent). (Notons que par hypothèse sur $g$ et d'après la proposition 2.1.3 la paire $(X', X'^\circ)$ est *torique*.)

$$\begin{array}{ccccc} T & \xleftarrow{g} & X & \xleftarrow{h} & Y \\ \uparrow & \square & \uparrow & \square & \uparrow \\ T' & \longleftarrow & X' & \xleftarrow{h'} & Y' \\ \uparrow & \square & \uparrow & \square & \uparrow \\ T'_i & \longleftarrow & X'_i & \xleftarrow{h'_i} & Y'_i \\ \downarrow & \square & \downarrow & \square \widetilde{h} & \downarrow \\ \widetilde{T} & \longleftarrow & \widetilde{X} & \xleftarrow{\widetilde{h}} & \widetilde{Y} \end{array}$$

Par stabilité par image inverse et théorème de changement de base propre, il en est de même des $\mathrm{R}^k h'^{\dagger}_{i\star}(\mathscr{F}_i')$. Enfin, la propriété d'une famille d'être constructible et modérée (resp. localement unipotente) étant Zariski locale, on a constructibilité et modération (resp. locale unipotence) des $\mathrm{R}^k \widetilde{h}^{\dagger}_{\star}(\widetilde{\mathscr{F}})$.

Quitte à remplacer $T$ par $\widetilde{T}$, on peut donc supposer la famille des $\mathrm{R}^k h^{\dagger}_{\star} \mathscr{F}$ constructible et modérée (resp. localement unipotente) *sur $X$*. (Notons que cela est plus fort que de supposer le résultat après altération de $X$.) La formule $\mathrm{R}f_{\star} = \mathrm{R}g_{\star} \circ \mathrm{R}h_{\star}$ (et ses variantes au-dessus des ouverts $T^{\dagger} = T[1/n]$) nous permet, via la suite spectrale associée, d'utiliser l'hypothèse de récurrence pour le triplet $(S, \rho, d-1)$ : la dimension de la fibre générique de $g$ est strictement inférieure à celle de $f$.

Ceci achève la démonstration du théorème 3.1.1, excepté le complément sur la majoration des fibres : on a montré que les énoncés pour les triplets $(S, \rho, d-1)$, $(S, \rho, -1, *)$, $(* \prec S, \rho, *)$ entraînent l'énoncé correspondant de modération ou locale unipotence pour le triplet $(S, \rho, d)$. On peut donc conclure par récurrence grâce au lemme du paragraphe 3.4.4. Initialisation : on rappelle que les énoncés sont trivialement vrais pour $S = \emptyset$ ou $\rho = -1$ et que l'énoncé pour un triplet $(S', \rho', -1)$ est également conséquence des $(* \prec S', \rho', *)$.

**3.5.4.** Notons que la *démonstration* du théorème 3.1.1 donne plus : pour tout changement de base *lisse* $R \to S$, la paire $(\beta : S' \twoheadrightarrow S, \mathfrak{S}) \times_S R$ — obtenue en prenant l'image inverse des données obtenues par la démonstration que nous avons suivie — convient pour le triplet $(f, \alpha, \mathfrak{X}) \times_S R$. Ceci résulte de la récurrence et du fait que si $(T, T^{\circ})$ est une paire régulière et $Q \to T$ lisse, la paire $(Q, Q^{\circ} := T^{\circ} \times_T Q)$ est encore régulière.

**3.5.5.** Mieux : la proposition 2.3.1 sous la forme la plus générale, dont la démonstration a été seulement *esquissée* en 2.3.6 dans le cas d'une base non nécessairement torique, permet d'établir le théorème d'uniformité propre dans le cas où $S$ est seulement nœthérien (et pas nécessairement quasi-excellent), et ce de façon compatible à *tout* changement de base (et pas seulement aux changements de base *lisses*). Ce raffinement sur les changements de base autorisés nous permet par descente des objets (morphismes, faisceaux, et stratifications) de supposer la base de type fini sur $\mathrm{Spec}(\mathbb{Z})$ donc en particulier excellente. Par normalisation (finie), on peut également la supposer normale. Les dévissages précédents s'appliquent sans changement et nous ramènent à la proposition susmentionnée pour les courbes nodales, sans utiliser le théorème de résolution par altération des singularités d'O. Gabber mais seulement le théorème de A. J. de Jong. (La récurrence 3.5.3 devient alors significativement plus simple : il s'agit seulement d'une récurrence sur la dimension relative dans le cas plurinodal.)

**3.6. Majoration des fibres.** Signalons rapidement en quoi la méthode suivie entraîne le complément, énoncé en 3.1.1, sur la « majoration » uniforme des fibres des images directes propres.

**3.6.1.** Réduction au cas $\alpha = \mathrm{Id}$ ($\leftrightarrow$ 3.3.1) : observer que si chaque $\mathrm{R}^i f_{j\star}^\dagger \mathscr{F}_j$ est $N_j$-extension ponctuelle de $\mathscr{F}_j$ relativement à $f_j^\dagger$, le faisceau $\mathrm{R}^k f_\star^\dagger \mathscr{F}$ est $N$-extension ponctuelle de $\mathscr{F}$ relativement à $f^\dagger$ pour $N \coloneqq (k+1) \times \max_{0 \leqslant j \leqslant 2d+1} N_j$.

**3.6.2.** Réduction au cas $\mathscr{F} = j_!^\dagger \mathscr{L}$ ($\leftrightarrow$ 3.3.3, 3.4.5) : d'une part, d'après 1.1.4, on a une borne *a priori* sur la longueur d'une filtration à gradués des prolongements par zéro de lisses et, d'autre part, la réduction au cas d'une strate ouverte est licite d'après 1.7.2 (ii).

**3.6.3.** Réduction au cas plurinodal ($\leftrightarrow$ 3.5.2) : utiliser (la récurrence et) la stabilité par extension, 1.7.2 (iii).

**3.6.4.** Réduction au cas nodal ($\leftrightarrow$ 3.5.3) : la propriété à démontrer est Zariski locale sur la source, stable par changement de base d'après 1.7.2 (v), et compatible à la factorisation $f = g \circ h$ au sens de 1.7.2 (iv).

**3.6.5.** Il reste à montrer que sous les hypothèses de la proposition 2.3.1, les faisceaux $\mathrm{R}^i h_\star(j_! \mathscr{L})$ ($0 \leqslant i \leqslant 2$) sont $N$-extensions ponctuelles de sous-quotients de $\mathscr{L}$ pour un entier $N$ ne dépendant que du morphisme $h : Y \to X$ et de l'immersion ouverte $j : Y^\circ \hookrightarrow Y$. Il suffit pour cela de vérifier que si $U$ est un ouvert d'une courbe projective lisse connexe $C$ sur un corps algébriquement clos, et $\mathscr{L}$ abélien lisse (de torsion inversible sur $C$) et modéré sur $U$, les groupes $\mathrm{H}^i_c(U, \mathscr{L})$ sont $N$-extensions ponctuelles de sous-quotients de $\mathscr{L}$ pour un entier $N$ ne dépendant que du *genre* de $C$ et du cardinal de $C - U$, qui varient l'un comme l'autre de façon constructible en famille. C'est évident pour $i = 0$, également vrai pour $i = 2$ par dualité et conséquence de la formule de Grothendieck-Ogg-Šafarevič [SGA 5 X, 7.1] pour $i = 1$. (Noter que l'on utilise cette formule dans le cas *modéré*, bien plus élémentaire ; cf. [OGG 1962, prop. 6, p. 199] ou [ŠAFAREVIČ 1961, th. 3, p. 337].)

# 4. Image directe par une immersion ouverte

**4.1. Énoncé et stratégie.** Le but de cette section est d'établir le théorème suivant.

**4.1.1. Théorème.** *Soient $X$ un schéma nœthérien quasi-excellent et $j : U \hookrightarrow X$ une immersion ouverte. La propriété d'une famille $\mathbb{F}$ de faisceaux sur $X$ d'être constructible et modérée (resp. constructible et localement unipotente) pour la topologie des altérations est préservée par le foncteur $\mathrm{R} j_\star j^\star$.*

*En d'autres termes, et plus précisément, pour tout entier $c \geqslant -1$, toute altération $\alpha : X_1 \to X$ et toute stratification $\mathfrak{X}_1$ de $X_1$, il existe une altération $\beta : X_2 \to X$ et une stratification $\mathfrak{X}_2$ de $X_2$ telles que pour tout entier $n \geqslant 1$ et tout faisceau $\mathscr{F}$ de $\mathbb{Z}/n\mathbb{Z}$-modules sur $X^\dagger = X[1/n]$ constructible et modéré (resp. constructible et localement unipotent) le long de $\mathfrak{X}_1^\dagger$ alors $\mathrm{R}^c j_\star^\dagger j^{\dagger\star} \mathscr{F}$ est constructible et modéré (resp. constructible et localement unipotent) le long de $\mathfrak{X}_2^\dagger$.*

*De plus :*

(a) *il existe un entier $N$, dépendant uniquement du quadruplet $(j, \alpha, \mathfrak{X}_1, c)$, tel que chaque $\mathrm{R}^c j_\star^\dagger j^{\dagger\star} \mathscr{F}$ comme ci-dessus soit $N$-extension ponctuelle de sous-quotients de $\mathscr{F}$ au sens de 1.7 ;*

(b) si $X$ est de type fini sur un schéma $S$, il existe un ouvert dense $S°$ de $S$ tel que la formation de chaque $\mathrm{R}^c j^\dagger_\star j^{\dagger\star}\mathscr{F}$ comme ci-dessus commute aux changements de base $S' \to S°$.

**4.1.2. Remarques.**

(i) Si $X$ est de dimension finie, il résulte de [STG XVIII$_\mathrm{B}$, 1.1] qu'il existe un entier $C$ tel que les $\mathrm{R}^c j^\dagger_\star j^{\dagger\star}\mathscr{F}$ comme ci-dessus soient nuls pour chaque $c \geqslant C$. On peut donc trouver des données $\beta$, $\mathfrak{X}_2$ et $N$ qui conviennent pour tout $c$. Il est vraisemblable qu'en adaptant la méthode de [STG XV], on puisse s'affranchir de cette hypothèse supplémentaire sur $X$. (Un point de départ possible est l'observation que [STG XV, 2.1.7] est valable plus généralement pour des coefficients lisses modérés.)

(ii) A priori, même si $\alpha = \mathrm{Id}$, on ne peut pas nécessairement prendre $\beta = \mathrm{Id}$.

**4.1.3. Stratégie.** On procède par récurrence sur l'entier $c$. La première des deux méthodes d'O. Gabber, exposée en [STG XIII, §2], s'adapte sans difficulté particulière si ce n'est qu'il nous faut établir un analogue du théorème de constructibilité générique [Th. finitude, 1.9 (i)] de P. Deligne sous une forme légèrement renforcée. (Voir la mise en garde 1.4.9, et 4.5 infra pour l'énoncé précis.)

**4.2. Dévissages.**

**4.2.1.** Pour alléger l'exposition, on ne fera en général pas mention des entiers $n \geqslant 1$, et on omet les foncteurs $-[1/n] = -^\dagger$ : on note $j$ pour $j^\dagger = j[1/n]$, etc., et $\Lambda$ pour $\mathbb{Z}/n\mathbb{Z}$.

Rappelons que l'on procède par récurrence sur le degré cohomologique $\leqslant c$, où $c$ est un entier $\geqslant -1$, en s'autorisant à changer $X$. Le cas $c = -1$ est tautologiquement vrai.

**4.2.2.** (*Réduction au cas $X_1 = X$.*) Cela résulte du théorème d'uniformité propre (3.1.1) et de la descente cohomologique, sous la forme $\mathrm{R}j_\star j^\star = \mathrm{R}\varepsilon_\star \mathrm{R}j_{\bullet\star} j_\bullet^\star$, où $\varepsilon : X_\bullet \to X$ est le 0-cosquelette de l'altération $\alpha : X_1 \to X$ et $j_\bullet : U_\bullet \hookrightarrow X_\bullet$ est le morphisme déduit de $j$ par changement de base $X_\bullet \to X$ (cf. 3.3.1),

Ce même argument nous permet de supposer que $X$ est un coproduit fini de schémas intègres et, finalement, que $X$ est intègre (cf. 3.3.2).

**4.2.3.** (*Réduction au cas où les faisceaux sont des prolongements par zéro de faisceaux lisses sur la strate ouverte $V \in \mathfrak{X}_1$ de $X$.*) Elle est entièrement semblable à la réduction effectuée en 3.3.3, une fois constaté que si $i : F \hookrightarrow X$ est une immersion fermée, on a $\mathrm{R}j_\star j^\star \circ i_\star = i_\star \circ \mathrm{R}(j \times_X i)_\star (j \times_X i)^\star$.

**4.2.4.** (*Réduction au cas où la strate ouverte $V$ est contenue dans $U$.*) D'après la réduction précédente, on peut supposer le faisceau $\mathscr{F}$ de la forme $h_!\mathscr{L}$ où $h$ est l'immersion ouverte de la strate ouverte $V$ dans $X$ et $\mathscr{L}$ est un faisceau lisse sur $V$ de $\Lambda$-modules. Notons $U°$ l'ouvert $V \cap U$.

$$\begin{array}{ccc} \mathscr{L} & & \\ V & \xrightarrow{h} & X \\ \uparrow & \square & \uparrow j \\ U^\circ & \xrightarrow{h^\circ} & U \\ \mathscr{L}^\circ := \mathscr{L}_{|U^\circ} & & \end{array}$$

Le complexe $\mathrm{R}j_\star j^\star(h_!\mathscr{L}) = \mathrm{R}j_\star(h_!^\circ \mathscr{L}^\circ)$ se réécrit $\mathrm{R}j_\star j^\star(j_! h_!^\circ \mathscr{L}^\circ)$. Quitte à raffiner la stratification, le sous-faisceau $j_! h_!^\circ \mathscr{L}^\circ$ de $h_!\mathscr{L}$ est constructible et modéré (resp. constructible et localement unipotent) le long de $\mathfrak{X}_1$ ; on peut donc supposer $V \subseteq U$, c'est-à-dire $V = U^\circ$.

**4.2.5**. *(Réduction au cas où $j^\star \mathscr{F}$ est lisse, c'est-à-dire $V = U$.)* Sous l'hypothèse précédente, le morphisme $h$ se factorise à travers $j$ selon le diagramme suivant :

$$\begin{array}{ccccc} \mathscr{L} & & h_!^\circ \mathscr{L} = j^\star \mathscr{F} & & \mathscr{F} = h_!\mathscr{L} \\ V = U^\circ & \xrightarrow{h^\circ} & U & \xrightarrow{j} & X. \\ & & \underset{h}{\searrow} & & \end{array}$$

On veut montrer que, sous nos hypothèses, les faisceaux $\mathrm{R}^c j_\star(h_!^\circ \mathscr{L})$ sont uniformément constructibles et modérés (resp. constructibles et localement unipotents) pour la topologie des altérations. Soient $F$ le fermé $U - U^\circ$ de $U$, $i : F \hookrightarrow U$ l'immersion fermée et considérons l'adhérence $\overline{F}$ de $F$ dans $X$, $j_F : F \hookrightarrow \overline{F}$ l'immersion ouverte correspondante et enfin $\overline{i}$ l'immersion fermée de $\overline{F}$ dans $X$.

$$\begin{array}{ccc} F & \xrightarrow{i} & U \\ j_F \downarrow & \overline{i} j_F = ji & \downarrow j \\ \overline{F} & \xrightarrow{\overline{i}} & X \end{array}$$

Le triangle distingué

$$\mathrm{R}j_\star(h_!^\circ \mathscr{L}) \longrightarrow \mathrm{R}h_\star \mathscr{L} \longrightarrow \mathrm{R}j_\star i_\star i^\star j^\star \mathrm{R}h_\star \mathscr{L} = \overline{i}_\star \mathrm{R}j_{F\star} j_F^\star (\overline{i}^\star \mathrm{R}h_\star \mathscr{L}) \xrightarrow{+1}$$

(déduit par application de $\mathrm{R}j_\star$ au triangle $h_!^\circ h^{\circ\star} \to \mathrm{Id} \to i_\star i^\star$ évalué en $\mathrm{R}h_\star^\circ \mathscr{L} = j^\star \mathrm{R}h_\star \mathscr{L}$) nous ramène à montrer l'uniformité (resp. la commutation générique aux changements de base) en degré $c$ du second terme et en degré $c-1$ du troisième. Par récurrence sur $c$, il suffit étudier le seul terme central : on peut donc dorénavant supposer que $j = h$, c'est-à-dire que $j^\star \mathscr{F}$ est égal au faisceau $\mathscr{L}$, lisse sur $U^\circ = U$.

**4.3. Cas d'un « bon » schéma.** Dans le cas particulier, souvent suffisant en pratique, où $X$ est de type fini sur un corps ou un schéma de Dedekind, la démonstration du théorème 4.1.1 est immédiate à partir de ce qui précède, et ne fait donc pas appel aux résultats génériques des paragraphes suivants. (On ne considère pas les compléments dans l'esquisse qui suit et renvoyons le lecteur à la démonstration donnée dans le cas général.)

En effet, d'après [DE JONG 1996, 4.1 et 8.2], il existe *pour un tel X* un hyperrecouvrement *propre* $X_\bullet \to X$ tel que la paire $(X_\bullet, U_\bullet)$ soit régulière, étage par étage. La descente cohomologique et le théorème d'uniformité dans le cas propre nous permettent donc de supposer la paire $(X, U)$ régulière. Or, il résulte de 2.4.1 que $\mathrm{R}j_\star \mathscr{L}$ est à cohomologie uniformément constructible et modérée (resp. constructible et localement unipotente) sur $X$. (Noter qu'il n'est pas nécessaire d'avoir recours à une altération du but.)

**4.4. Cas général : réduction à l'uniformité générique.** Traitons maintenant le cas modéré, sans supposer le schéma $X$ de type fini sur un corps ou un schéma de Dedekind. Le cas localement unipotent est semblable et laissé au lecteur ; la majoration des fibres sera brièvement traitée ci-après (4.5.7). Quant à la propreté cohomologique générique, elle résulte immédiatement du théorème de changement de base propre et de 4.5.3 (ii).

Une récurrence nœthérienne sur (les fermés de) $X$, le lemme de Gabber 1.6.2 et l'observation 1.4.10 appliquée aux faisceaux $\mathscr{K} := \mathrm{R}^c j_\star \mathscr{L}$ (pour $\mathscr{L}$ lisse sur $U$ comme précédemment) nous ramènent à montrer que si $i : F \hookrightarrow X$ est un fermé contenu dans $X \smallsetminus U$, il existe une immersion ouverte $k : F^\circ \hookrightarrow F$ d'image dense, une altération $a : F' \to F$ et enfin une stratification $\mathfrak{F}$ de $F$ telles que les faisceaux $k_! k^\star(\mathscr{K}_{|F})$ soient constructibles et modérés le long de $\mathfrak{F}$. D'après le théorème de résolution des singularités d'O. Gabber ([STG VII, 1.1]), il existe un hyperrecouvrement $\varepsilon_X : X_\bullet \to X$ pour la topologie $h$ tel qu'à chaque étage $n$, le schéma $X_n$ soit régulier, $F_n := F \times_X X_n$ un diviseur à croisements normaux strict et $U_n := U \times_X X_n$ le complémentaire d'un tel diviseur. D'après la formule d'hyperchangement de base ([STG XII$_\mathrm{A}$, 2.2.5],[STG XII$_\mathrm{B}$]) d'O. Gabber, le morphisme

$$i^\star \mathrm{R} j_\star \mathscr{L} \xrightarrow{\sim} \mathrm{R}\varepsilon_{F\star}(i_\bullet^\star \mathrm{R} j_{\bullet\star} \mathscr{L}_\bullet)$$

est un isomorphisme. Les morphismes $\varepsilon_{F,n} : F_n \to F$ n'étant pas nécessairement propres, on ne peut appliquer directement le théorème d'uniformité propre. Cependant, la proposition 4.5.1 ci-dessous. appliquée pour $n$ variable aux morphismes $f := \varepsilon_{F,n}$ (dans le cas particulier $S = Y$), nous ramène à démontrer le fait suivant, lorsque $X' = X_n$, $F' = F_n$, et $U' = U_n$.

**Lemme.** *Soient*

$$\begin{array}{ccccccc}
U' & \xhookrightarrow{j'} & X' & \xleftarrow{i'} & F' & \xhookdashrightarrow{k'} & \overline{F}' \\
\downarrow & & \downarrow & & \downarrow & \swarrow_{propre} & \\
U & \xhookrightarrow{j} & X & \xleftarrow{i} & F & &
\end{array}$$

*un diagramme commutatif à carrés cartésiens et traits pleins, où :*

— *$X'$ est régulier, séparé de type fini sur $X$ ;*
— *$F'$ est un diviseur à croisements normaux strict ;*
— *$U'$ est le complémentaire d'un tel diviseur.*

*Alors, il existe une $F$-compactification $k' : F' \hookrightarrow \overline{F}'$ de $F'$ telle que pour tout faisceau lisse $\mathscr{L}$ de $\mathbb{Z}/n\mathbb{Z}$-modules sur $U$ tel que $j_!\mathscr{L}$ soit modéré (resp. localement unipotents) et $n \geq 1$ inversible sur $X$, les faisceaux $k'_! i'^{\star} \mathrm{R}^{c'} j'_\star \mathscr{L}$, pour $c' \leq c$, sont modérés et constructibles (resp. localement unipotents) le long d'une stratification de $\overline{F}'$.*

*Démonstration du lemme.* Il suffit de considérer une compactification de $X' \to X$ — à l'infini de laquelle le faisceau $j'_!\mathscr{L}$ est automatiquement modéré car $j_!\mathscr{L}$ l'est sur $X$ — et la compactification induite de $F'$ par changement de base $F \hookrightarrow X$. La conclusion résulte de 2.4.1 (ii). □

4.5. **Énoncés génériques.** Le but de ce paragraphe est d'établir la proposition suivante et d'achever ainsi la démonstration de 4.1.1 (excepté le (a), traité ultérieurement). Pour alléger les énoncés, on fixe l'entier $c$ et le schéma nœthérien quasi-excellent $S$.

4.5.1. **Proposition.** *Soient $f : X \to Y$ un morphisme séparé de type fini entre $S$-schémas de type fini et $j : X \hookrightarrow \overline{X}$ une $Y$-compactification de $f$. Alors, si $\mathbb{F}$ est une famille sur $X$ telle que $j_!\mathbb{F}$ soit constructible modérée (resp. localement unipotente) sur $\overline{X}$, il existe une immersion ouverte dominante $k_S : S^\circ \hookrightarrow S$ telle que $k_{Y!}\mathrm{R}^c f^\circ_\star \mathbb{F}$ soit constructible modérée (resp. localement unipotente) pour la topologie des altérations sur $S$, où $f^\circ := f \times_S S^\circ$ et $k_Y := k \times_S Y$.*

*Plus précisément, pour toute stratification $\mathfrak{X}$ de $\overline{X}$, il existe une altération $\beta : Y' \to Y$ et une stratification $\mathfrak{Y}$ de $Y'$ telles que pour tout entier $n \geq 1$ et tout faisceau $\mathscr{F}$ de $\mathbb{Z}/n\mathbb{Z}$-modules sur $X^\dagger = X[1/n]$ tel que $j^\dagger_! \mathscr{F}$ soit constructible et modéré (resp. constructible et localement unipotent) le long de $\mathfrak{X}^\dagger$ alors chaque $\beta^{\dagger\star} k^\dagger_{Y!} \mathrm{R}^c f^{\circ\,\dagger}_\star \mathscr{F}$ est constructible et modéré (resp. constructible et localement unipotent) le long de $\mathfrak{Y}^\dagger$. (On note abusivement $\mathscr{F}$ pour sa restriction à $X^{\circ\dagger}$.)*

*De plus, la formation de $\mathrm{R}^c f^{\circ\,\dagger}_\star \mathscr{F}$ commute aux changements de base $S' \to S^\circ$.*

4.5.2. **Remarques.**
(a) L'énoncé ne dit pas que si $\mathbb{F}$ est une famille constructible modérée (resp. localement unipotente) sur $X$, il existe une immersion ouverte dominante $Y^\circ \hookrightarrow Y$ provenant de $S$ telle que $\mathrm{R}f^\circ_\star \mathbb{F}$ soit constructible modérée (resp. localement unipotente) sur $S^\circ$ : on fait une hypothèse à l'infini (relativement à $Y$). (Voir 0.3 pour un contre-exemple à un tel énoncé.)
(b) Insistons sur le fait que cet énoncé est un corollaire du théorème 4.1.1, d'après lequel on peut prendre $S^\circ = S$ (c'est-à-dire $k = \mathrm{Id}$).

En particulier :

4.5.3. **Corollaire.** *Soient $X$ un $S$-schéma de type fini et $j : U \hookrightarrow X$ une immersion ouverte. Pour toute stratification $\mathfrak{X}_1$ de $X$, il existe une altération $\beta : X_2 \to X$, une stratification $\mathfrak{X}_2$ de $X_2$ et une immersion ouverte dominante $k_S : S^\circ \hookrightarrow S$ telles que pour tout entier $n \geq 1$ et tout faisceau $\mathscr{F}$ de $\mathbb{Z}/n\mathbb{Z}$-modules sur $X^\dagger = X[1/n]$ constructible et modéré (resp. constructible et localement unipotent) le long de $\mathfrak{X}^\dagger_1$ alors :*
  (i) *chaque $k^\dagger_{X!} k^{\dagger\star}_X \mathrm{R}^c j^\dagger_\star j^{\dagger\star} \mathscr{F}$ est constructible et modéré (resp. constructible et localement unipotent) le long de $\mathfrak{X}^\dagger_2$, où $k_X$ est l'immersion ouverte $k_S \times_S X$ ;*
  (ii) *la formation de chaque $\mathrm{R}^c j^\dagger_\star j^{\dagger\star} \mathscr{F}$ commute aux changements de base $S' \to S^\circ$.*

**4.5.4.** Montrons que, réciproquement, le cas particulier d'une immersion ouverte entraîne la proposition 4.5.1, dont nous reprenons les notations. Soit $k_S : S° \hookrightarrow S$ une immersion ouverte comme en 4.5.3 associée à l'immersion ouverte $j : X \hookrightarrow \overline{X}$. Par changement de base, on obtient la moitié droite du diagramme commutatif ci-dessous.

$$\begin{array}{c}
\xymatrix{
& & k_X & & \\
X \ar[r]^{j} \ar[d]_{f} & \overline{X} \ar[r]^{k_{\overline{X}}} \ar[d]^{\overline{f}} & \overline{X}° \ar[r]^{j°} \ar[d]^{\overline{f}°} & X° \ar[d]^{f°} \\
& Y \ar[r]^{k_Y} & Y° & \\
& S \ar[r]^{k_S} & S° &
}
\end{array}$$

On s'intéresse au foncteur $k_{Y!}f°_\star k_X^\star$ (où l'on omet les R). Compte tenu des isomorphismes $f°_\star = \overline{f}°_\star j°_\star$ (transitivité), $k_{Y!}\overline{f}°_\star = \overline{f}_\star k_{\overline{X}!}$ (propreté de $\overline{f}$) et $k_{\overline{X}}^\star j_\star = j°_\star k_X^\star$ (changement de base trivial), on a donc

$$k_{Y!}f°_\star k_X^\star = \overline{f}_\star \circ k_{\overline{X}!}k_{\overline{X}}^\star \circ j_\star.$$

La conclusion résulte aussitôt de l'hypothèse sur $k_{\overline{X}!}k_{\overline{X}}^\star j_\star$ et du théorème d'uniformité propre (3.1.1).

La commutation aux changements de base de $f°_\star = \overline{f}°_\star j°_\star$ résulte du théorème de changement de base propre pour $\overline{f}°$ et de la propreté cohomologique de $j°$ (relativement aux faisceaux considérés).

La suite de ce 4.5 est donc consacrée à la démonstration du corollaire 4.5.3.

**4.5.5.** *Dévissages.*
*Mutatis mutandis*, les réductions suivant 4.2.2 s'appliquent dans ce cadre générique. Pour la dernière, au cas du prolongement par zéro d'un faisceau lisse sur la strate dense, on utilise l'hypothèse de récurrence selon laquelle l'énoncé 4.1.1 est vrai pour l'entier $c - 1$.

(On aurait pu faire intervenir dans l'énoncé du corollaire une altération $\alpha : X_1 \to X$ non nécessairement identité ; la réduction au cas $\alpha = \mathrm{Id}$ aurait alors été possible par commutation de $k_! k^\star$ aux images directes propres.)

**4.5.6.**
**4.5.6.1.** Démontrons 4.5.3, dans le cas particulier — suffisant d'après ce qui précède — où les faisceaux test $\mathcal{F}$ sont de la forme $j_! \mathcal{L}$, avec $\mathcal{L}$ lisse sur $U$ (ou plus précisément $U^\dagger$). Nous allons faire un usage *relatif* de [DE JONG 1996] nous permettant — à l'instar de [ILLUSIE 2010, §3] et [ORGOGOZO 2003, §2] — de supposer que $U$ est génériquement le complémentaire d'un diviseur à croisements normaux (relatif) dans un $S$-schéma lisse, comme nous allons maintenant l'expliquer. Commençons par observer que le problème étant local sur $X$ et $S$, on peut supposer le morphisme de type fini $X \to S$ séparé, et $S = \mathrm{Spec}(A)$ intègre, de point générique $\eta = \mathrm{Spec}(K)$. Comme expliqué dans les deux références ci-dessus, on contourne la difficulté due au fait que $K$ n'est pas parfait par passage à la limite. Soit la perfection « générique » $T$ de $S$ : si $\mathrm{car.}(K) = p > 0$, c'est le spectre de la $\mathbb{F}_p$-algèbre $A^{\mathrm{parf}} = \mathrm{colim}_{a \mapsto a^p} A$ parfaite[10], ainsi que ses localisés ; sinon, $T \coloneqq S$. Les constructions géométriques (en nombre fini) suivantes, *a priori* faites sur $T$ pour pouvoir appliquer le théorème de résolution des singularités de A. J. de Jong ([DE JONG 1996, 4.1-2]) au-dessus de son point générique parfait, se descendent à un schéma $S' = \mathrm{Spec}(B)$ universellement homéomorphe à $S$. Un tel morphisme est innocent ; on peut donc supposer abusivement que $S = T$ (et en particulier $K = K^{\mathrm{parf}}$) pour simplifier les notations.

**4.5.6.2.** La technique de [ILLUSIE 2010] et [ORGOGOZO 2003] conduit pour chaque entier $N$ à un « bon » hyperrecouvrement propre (tronqué) $X^\circ_{\bullet \leqslant N} \to X^\circ \coloneqq X \times_S S^\circ$ au-dessus d'un ouvert dense $S^\circ$ de $S$. « Bon » signifie ici : lisse sur $S^\circ$, étage par étage, et tel que les $U^\circ_n \coloneqq U \times_X X^\circ_n$ soient complémentaires d'un diviseur à croisements normaux strict relativement à $S^\circ$. On souhaite ici un hyperrecouvrement propre $X_{\bullet \leqslant N} \to X$ induisant sur un ouvert dense $S^\circ$ un hyperrecouvrement du type précédent. Commençons par énoncer un lemme, conséquence immédiate du théorème [DE JONG 1996, 4.1] susmentionné.

**Lemme.** *Soit le diagramme cartésien à traits pleins ci-dessous.*

$$\begin{array}{ccccc}
W_\eta & \dashrightarrow & Z_\eta & \dashrightarrow & Z \\
\downarrow & & \downarrow & & \downarrow \\
V_\eta & \hookrightarrow & Y_\eta & \longrightarrow & Y \\
\downarrow & & \downarrow & & \downarrow \\
U_\eta & \hookrightarrow & X_\eta & \longrightarrow & X & \hookleftarrow & U \\
& & \downarrow & & \downarrow \\
& & \eta & \longrightarrow & S
\end{array}$$

*Il existe un morphisme propre et surjectif $Z \to Y$ induisant les carrés cartésiens supérieurs tel que $Z_\eta$ soit lisse sur $\eta$ et que l'ouvert $W_\eta \subseteq Z_\eta$ en soit le complémentaire d'un diviseur à croisements normaux strict.*

(Prendre pour $Z$ une $Y$-compactification d'une résolution de $Y_\eta$ adaptée à l'ouvert $V_\eta$.)

Les techniques classiques de construction d'hyperrecouvrements — voir aussi [GILLET et SOULÉ 2009, §2.4] pour un problème semblable — nous permettent de déduire de ce lemme une construction étage par étage d'un hyperrecouvrement *propre* tronqué $X_{\bullet \leqslant N} \to X$ tel que les paires $(X_n, U_n)$ soient génériquement régulières : $X_{\bullet \leqslant n}$ étant construit, on pose $Y \coloneqq \mathrm{cosq}(X_{\bullet \leqslant n})_{n+1}$ et $X_{n+1}$ est alors le $Y$-schéma propre $Z$ fourni par le lemme. La lissité au-dessus de $\eta$ des schémas $X_n$ et des branches des diviseurs $X_n \smallsetminus U_n$ ainsi construits pour $n \leqslant N$ s'étend à un ouvert dense $S^\circ$ de $S$.

---

[10]Mais pas nécessairement nœthérienne.

**4.5.6.3.** En omettant les troncations, on a donc un diagramme commutatif à carrés cartésiens

$$
\begin{array}{ccccccc}
U_\bullet^\circ & \xhookrightarrow{j_\bullet^\circ} & X_\bullet^\circ & \xrightarrow{k_{X_\bullet}} & X_\bullet & \xhookleftarrow{j_\bullet} & U_\bullet \\
\downarrow & & \downarrow{\varepsilon^\circ} & & \downarrow{\varepsilon} & & \downarrow \\
U^\circ & \hookrightarrow & X^\circ & \xrightarrow{k_X} & X & \xhookleftarrow{j} & U \\
& \searrow & \downarrow & & \downarrow & \swarrow & \\
& & S^\circ & \hookrightarrow & S & &
\end{array}
$$

Par descente cohomologique, on peut réécrire le foncteur étudié :

$$k_{X!}k_X^\star \circ j_\star j^\star = \varepsilon_\star \circ k_{X_\bullet!}k_{X_\bullet}^\star \circ j_{\bullet\star}j_\bullet^\star \circ \varepsilon^\star = \varepsilon_\star \circ k_{X_\bullet!}j_{\bullet\star}^\circ \circ (U_\bullet \to X)^\star.$$

Or, il résulte de 2.4.1 (ii) que pour $\mathscr{L}$ lisse sur $U$ tel que $j_!\mathscr{L}$ soit modéré (resp. localement unipotent), les faisceaux $k_{X_n!}\mathrm{R}^c j_\star^\circ \mathscr{L}^\circ$ le sont également. L'énoncé 4.5.3 (i) résulte alors du théorème d'uniformité propre pour les $\varepsilon_n$, $n \leqslant c < N$. Enfin, la propreté cohomologique des $j_n^\circ$ relativement à $S^\circ$ ayant été établie en 2.4.1 (cf. 2.4.3), l'énoncé 4.5.3 (ii) résulte du théorème de changement de base propre (pour les $\varepsilon_n^\circ$).

**4.5.7**. *Majoration des fibres.* Pour établir 4.1.1 (a), la démonstration précédente nous ramène sans difficulté à la proposition 2.4.6, déjà établie et à l'énoncé dans le cas propre. Dans le cas « absolu », passant par les énoncés génériques, on procède à nouveau par récurrence sur les fermés de $X \smallsetminus U$ ; la mise en garde 1.4.9 sur la nécessité de considérer des prolongements par zéro n'a par contre trivialement pas lieu d'être.

## 5. Six opérations

**5.1. Image directe par un morphisme non nécessairement propre.** Le théorème suivant résulte immédiatement, par composition, des énoncés 3.1.1 (cas d'un morphisme propre) et 4.1.1 (cas d'une immersion ouverte).

**Théorème.** *Soient $S$ un schéma nœthérien quasi-excellent, $f : X \to S$ un morphisme de type fini, $j : X \hookrightarrow \overline{X}$ une $S$-compactification de $X$ et $\mathbb{F}$ une famille de faisceaux sur $X$ telle que la famille $j_!\mathbb{F}$ soit constructible et modérée (resp. constructible et localement unipotente) pour la topologie des altérations. Alors, la famille des complexes $\mathrm{R}f_\star\mathbb{F}$ sur $S$ est constructible et modérée (resp. constructible et localement unipotente) pour la topologie des altérations.*

$$\begin{array}{c}
\overline{\mathfrak{X}}\ \overline{X}' \\
\downarrow \alpha \\
X \xhookrightarrow{j} \overline{X} \\
f \downarrow \quad \swarrow \overline{f} \\
S \xleftarrow{\exists \beta} S' \\
\exists \mathfrak{S}
\end{array}$$

*En d'autres termes, et plus précisément, pour tout entier $c \geqslant 0$, toute altération $\alpha : \overline{X}' \to \overline{X}$ et toute stratification $\overline{\mathfrak{X}}$ de $\overline{X}'$, il existe une altération $\beta : S' \to S$ et une stratification $\mathfrak{S}$ de $S'$ telles que pour tout entier $n \geqslant 1$ et tout faisceau $\mathscr{F}$ de $\mathbb{Z}/n\mathbb{Z}$-modules sur $X^\dagger = X[1/n]$ tel que $\alpha^{\dagger\star} j_!^\dagger \mathscr{F}$ soit constructible et modérée (resp. constructible et localement unipotent) le long de $\overline{\mathfrak{X}}^\dagger$ alors chaque $\beta^{\dagger\star} \mathrm{R}^c f_\star^\dagger \mathscr{F}$ est constructible et modéré (resp. constructible et localement unipotent) le long de $\mathfrak{S}^\dagger$.*

*De plus, il existe un entier $N$ et un ouvert dense $S^\circ \subseteq S$, dépendant uniquement du quadruplet $(f, \alpha, \overline{\mathfrak{X}}, c))$, tels que chaque $\mathrm{R}^c f_\star^\dagger \mathscr{F}$ comme ci-dessus*

(i) *soit $N$-extension ponctuelle de sous-quotients de $\mathscr{F}$ ;*
(ii) *commute aux changements de base $T \to S^\circ$.*

**Remarques.**
(a) L'hypothèse d'excellence est ici essentielle.
(b) Pour la dépendance en $c$, voir 4.1.2 (i).
(c) On dit pas que $\mathrm{R}^c f_\star^\dagger \mathscr{F}$ est $N$-extension ponctuelle de $\mathscr{F}$ relativement à $f$ (au sens de 1.7) : c'est trivialement faux lorsque $f = j \neq \mathrm{Id}$ par exemple.

### 5.2. Autres opérations.

**5.2.1.** Nous avons établi un résultat de stabilité par $f_\star$, en faisant une hypothèse — dont on ne peut totalement s'affranchir — sur le comportement des faisceaux à l'infini d'une même compactification.

**5.2.2.** Cette même hypothèse rend tautologique la stabilité par $f_!$ à partir du théorème d'uniformité propre.

Considérons maintenant brièvement les 4 des 6-opérations restantes.

**5.2.3.** La stabilité par $f^\star$ est tautologique.

**5.2.4.** Comme expliqué en [Beilinson, Bernstein et Deligne 1982, p. 151] ou [Th. finitude, 1.5], la stabilité par $f^!$ se ramène à la stabilité par $\mathrm{R} j_\star j^\star$, déjà établie : le théorème ci-dessous résulte des faits suivants :

— la propriété à démontrer est locale pour la topologie de Zariski (en haut et en bas) ;
— un morphisme affine est lissifiable ;
— si $S' \to S$ est une altération, il en est de même de $\mathbb{A}^d_{S'} \to \mathbb{A}^d_S$ ;
— dualité de Poincaré et innocuité de la torsion à la Tate ;
— existence d'un triangle distingué $i_\star i^! \to \mathrm{Id} \to j_\star j^\star \xrightarrow{+1}$.

**Théorème.** *Soient $S$ un schéma nœthérien quasi-excellent et $f : X \to S$ un morphisme de type fini. La propriété d'une famille $\mathbb{F}$ de faisceaux (sur $S$) d'être constructible et modérée (resp. constructible et localement unipotente) pour la topologie des altérations est préservée par le foncteur $f^!$.*

*En d'autres termes, et plus précisément, pour tout entier $c \geqslant 0$, toute altération $\alpha : S' \to S$ et toute stratification $\mathfrak{S}$ de $S'$, il existe une altération $\beta : X' \to X$ et une stratification $\mathfrak{X}$ de $X'$ telles que pour tout entier $n \geqslant 1$ et tout faisceau $\mathscr{F}$ de $\mathbb{Z}/n\mathbb{Z}$-modules sur $S^\dagger = S[1/n]$ tel que $\alpha^{\dagger\star}\mathscr{F}$ soit constructible et modérée (resp. constructible et localement unipotent) le long de $\mathfrak{S}^\dagger$, chaque $\beta^{\dagger\star}\mathrm{R}^c f^{\dagger!}\mathscr{F}$ est constructible et modéré (resp. constructible et localement unipotent) le long de $\mathfrak{X}^\dagger$.*

$$\begin{array}{ccc}
 & \exists\mathfrak{X} \\
X & \xleftarrow{\exists\beta} & X' \\
& \downarrow f & \\
\mathscr{F}\ S & \xleftarrow{\alpha} & S' \\
& & \mathfrak{S}
\end{array}$$

**5.2.5.** *Stabilité par $\otimes^{\mathbb{L}}$.* On souhaite montrer que si $(X, \mathfrak{X})$ est un schéma nœthérien stratifié alors pour tout $n \geqslant 1$ et toute paire $\mathscr{A}, \mathscr{B} \in \mathrm{Ob}\, D_{\mathfrak{X}}^{-\mathrm{mod}}(X[1/n], \mathbb{Z}/n\mathbb{Z})$ (variante évidente de la notation introduite en 3.1.3) le produit tensoriel $\mathscr{A} \otimes^{\mathbb{L}}_{\mathbb{Z}/n\mathbb{Z}} \mathscr{B}$ appartient également à $D_{\mathfrak{X}}^{-\mathrm{mod}}(X[1/n], \mathbb{Z}/n\mathbb{Z})^{\tiny\textcircled{1}}$ et de même pour les complexes à cohomologie localement unipotente le long de $\mathfrak{X}$.

Par dévissage des complexes, on peut supposer que $\mathscr{A}$ et $\mathscr{B}$ sont des faisceaux placés en degré nul. Par dévissage (1.1.4) des faisceaux constructibles le long de $\mathfrak{X}$ — ou plutôt $\mathfrak{X}^\dagger = \mathfrak{X}[1/n]$ mais nous ignorons cette subtilité pour simplifier les notations —, de l'entier $n$ (voir par exemple [Th. finitude, A, 2.2 b)] ou [STG XVII, 7.4.6] pour un raffinement), et grâce à la formule de projection triviale (pour une immersion ouverte ou fermée), on se ramène à l'étude des produits tensoriels $j_!(\mathscr{L} \otimes^{\mathbb{L}}_{\mathbb{Z}/\ell^r\mathbb{Z}} \mathscr{M})$, où $\mathscr{L}, \mathscr{M}$ sont des faisceaux lisses de $\mathbb{F}_\ell$-espaces vectoriels sur une même strate ouverte dense $U \subseteq X = X[1/\ell]$. (Cette réduction n'utilise pas le théorème d'uniformité propre pour une immersion fermée — c'est-à-dire essentiellement le lemme de Gabber — car il n'est pas nécessaire d'avoir recours à une altération.)

Comme d'autre part $\mathscr{L} \otimes^{\mathbb{L}}_{\mathbb{Z}/\ell^r\mathbb{Z}} \mathscr{M} = (\mathscr{L} \otimes^{\mathbb{L}}_{\mathbb{Z}/\ell^r\mathbb{Z}} \mathbb{F}_\ell) \otimes^{\mathbb{L}}_{\mathbb{F}_\ell} \mathscr{M}$ est représenté (pour $r > 1$) par le complexe (concentré en degrés $\leqslant 0$)

$$[\cdots \to \mathscr{L} \otimes_{\mathbb{F}_\ell} \mathscr{M} \xrightarrow{0} \mathscr{L} \otimes_{\mathbb{F}_\ell} \mathscr{M} \xrightarrow{0} \mathscr{L} \otimes_{\mathbb{F}_\ell} \mathscr{M} \xrightarrow{0} \mathscr{L} \otimes_{\mathbb{F}_\ell} \mathscr{M}],$$

il suffit d'observer que si $j_!\mathscr{L}$ et $j_!\mathscr{M}$ sont modérés (resp. localement unipotents le long d'une stratification $\mathfrak{X}$), il en est de même de $j_!(\mathscr{L} \otimes_{\mathbb{F}_\ell} \mathscr{M})$.

**5.2.6.** *Stabilité par* **RHom**. On suppose $X$ de dimension finie, pour éviter d'avoir *a priori* besoin de tronquer (cf. 4.1.2). Soient $\alpha : X_1 \to X$ une altération et $\mathfrak{X}_1$ une stratification de $X_1$. On veut montrer qu'il existe une altération $\beta : X_2 \to$ $X$ et une stratification $\mathfrak{X}_2$ de $X_2$ telles que pour chaque $n$, le foncteur **RHom** : $D(X[1/n], \mathbb{Z}/n\mathbb{Z})^{\mathrm{op}} \times D(X[1/n], \mathbb{Z}/n\mathbb{Z}) \to D(X[1/n], \mathbb{Z}/n\mathbb{Z})$ induise un foncteur

$$D^{-\mathrm{mod}}_{\alpha, \mathfrak{X}_1}(X[1/n], \mathbb{Z}/n\mathbb{Z})^{\mathrm{op}} \times D^{+\mathrm{mod}}_{\alpha, \mathfrak{X}_1}(X[1/n], \mathbb{Z}/n\mathbb{Z}) \to D^{+\mathrm{mod}}_{\beta, \mathfrak{X}_2}(X[1/n], \mathbb{Z}/n\mathbb{Z}).$$

---

[1] On ne peut, bien sûr, se limiter aux catégories dérivées *bornées* : si $r > 1$, le complexe $\mathbb{F}_\ell \otimes^{\mathbb{L}}_{\mathbb{Z}/\ell^r\mathbb{Z}} \mathbb{F}_\ell$ est dans $D^-(\mathbb{Z}/\ell^r\mathbb{Z})$ mais pas $D^b(\mathbb{Z}/\ell^r\mathbb{Z})$.

On procède à nouveau comme dans [BEILINSON, BERNSTEIN et DELIGNE 1982, 6.1.2, p. 151]. (Voir aussi [Th. finitude, A, 1.6], [SGA 5 I, §3.3] ou [STG XVII, §7.6].) La formule

$$\mathrm{R}f_\star \mathbf{RHom}(A, f^! B) = \mathbf{RHom}(f_! A, B)$$

pour une immersion $f$ nous ramène par dévissages et stabilité par $\mathrm{R}j_\star$ (4.1.1) — nécessitant peut-être de grossir $X_1$ et raffiner $\mathfrak{X}_1$ — à montrer que si $\mathscr{L}$ est un faisceau lisse sur un ouvert dense $j : U \hookrightarrow X$ tel que $j_!\mathscr{L}$ soit modéré (resp. localement unipotent) le long de $\mathfrak{X}_1$, et $\mathscr{G}$ est également constructible et modéré (resp. localement unipotent) le long de $\mathfrak{X}_1$, alors il en est de même des faisceaux $j_! \mathbf{Ext}^i(\mathscr{L}, j^\star \mathscr{G})$ ($i \geqslant 0$). Or, sous l'hypothèse faite sur $\mathscr{L}$, il existe un carré cartésien

$$\begin{array}{ccc} U' & \xhookrightarrow{j'} & X' \\ \downarrow & & \downarrow \\ U & \xhookrightarrow{j} & X \end{array}$$

avec $X' \to X$ fini surjectif, $U' \to U$ fini étale tel que $\mathscr{L}_{|U'}$ soit trivial et si $\mathscr{H}$ est un faisceau sur $U_1 \coloneqq U \times_X X_1$, la constructibilité et modération (resp. locale unipotence) de $j'_{1!}\mathscr{H}'_1$ (sur $X'_1 \coloneqq X_1 \times_X X'$) le long de $\mathfrak{X}'_1$ entraîne la constructibilité et modération (resp. locale unipotence) de $j_{1!}\mathscr{H}_1$ le long de $\mathfrak{X}_1$. (Dans le cas unipotent, on utilise le fait que si un groupe $G$ agit sur un $\mathbb{F}_\ell$-module $M$ de sorte que la restriction à $H \triangleleft G$ soit unipotente alors l'action de $G$ est également unipotente si le quotient $G/H$ est un $\ell$-groupe.[12]) Écrivant $\mathscr{L}$ comme un quotient de $(U' \to U)_\star (U' \to U)^\star \mathscr{L}$ et utilisant la formule d'adjonction rappelée ci-dessus (pour $f$ le morphisme étale $U' \to U$), ceci nous ramène au cas particulier où $\mathscr{L}$ est même trivial, puis par résolution projective au cas $\mathscr{L} = \Lambda$. Pour conclure, il ne reste plus qu'à utiliser le fait que $\mathbf{RHom}_\Lambda(\Lambda, j^\star \mathscr{G}) = j^\star \mathscr{G}$.

## 6. APPLICATIONS

**6.1. Représentations $\ell$-adiques.** Dans [BÖCKLE, GAJDA et PETERSEN 2015], les auteurs démontrent une conjecture de *presque indépendance* des représentations $\ell$-adiques associées à un schéma algébrique, due à J.-P. Serre en caractéristique nulle ([SERRE 2013, §3.2]). Un des ingrédients de leur démonstration est le fait suivant ([BÖCKLE, GAJDA et PETERSEN 2015, 6.3, 7.4]), que les auteurs démontrent par réduction au cas des courbes (à la Wiesend) et en utilisant [DELIGNE 1973, 9.8] :

**Proposition.** *Soit $k$ un corps parfait, $K/k$ une extension de type finie et $X_K$ un $K$-schéma algébrique. Il existe une extension finie étale $L/K$, un $k$-schéma propre et lisse $T$ de corps des fractions $L$ et un diviseur à croisements normaux $D$ dans $T$ tels que* pour chaque nombre premier $\ell$ inversible sur $k$,

---

[12] En effet, $\mathrm{Fix}(H \curvearrowright M) \neq 0$ est stable par $G$ donc muni d'une action du $\ell$-groupe fini $G/H$. Une telle action a un point fixe non trivial.

(i) *l'action du groupe de Galois $G_L$ sur $\mathrm{H}^\star(X_{\overline{K}}, \mathbb{Q}_\ell)$ se factorise à travers le groupe fondamental* modéré $\pi_1^{\mathrm{mod}}(T, D)$ ;
(ii) *pour chaque point maximal $d$ de $D$, l'action sur $\mathrm{H}^\star(X_{\overline{K}}, \mathbb{Q}_\ell)$ du groupe d'inertie correspondant soit* unipotente.

C'est également un corollaire immédiat de nos résultats.

*Démonstration.* Par déploiement, il existe un $k$-schéma de type fini $U$ de corps des fractions $K$, et un morphisme de type fini $X \to U$ de fibré générique $X_K$. Soit $S$ une $k$-compactification de $U$ et notons $f : X \to S$ le morphisme composé. D'après 5.1, il existe une altération $\beta : T \to S$ et une stratification $\mathfrak{T}$ de $T$ telle que les faisceaux $\beta^*\mathrm{R}^i f_\star \mathbb{Q}_\ell$, pour $\ell \neq \mathrm{car.}(k)$ et $i \in \mathbb{N}$, soient constructibles, modérés et localement unipotents le long de $\mathfrak{T}$. Quitte à rétrécir $U$, on peut supposer ces faisceaux lisses sur cet ouvert et que $V \coloneqq \beta^{-1}(W)$ est la strate ouverte dense de $T$. Comme remarqué en 1.6.6, on peut supposer l'altération $\beta$ génériquement étale. D'après [DE JONG 1996, 4.1], on peut également supposer le morphisme $T \to \mathrm{Spec}(k)$ lisse et l'ouvert $V$ être le complémentaire d'un diviseur à croisements normaux $D$. □

**Remarque.** Le résultat de P. Deligne ([BERTHELOT 1997, 6.3.2]) d'indépendance de $\ell$ dans le théorème de monodromie de Grothendieck est également un corollaire immédiat de 5.1, du moins pour un trait excellent. Pour des précisions numériques dans cette direction, voir [UMEZAKI 2016].

**6.2. Lien avec la torsion.** Soient $k$ un corps algébriquement clos d'exposant caractéristique $p \geqslant 1$ et $X$ un schéma projectif et lisse sur $k$. Pour chaque nombre premier $\ell$ différent de $p$ et chaque entier $i \geqslant 0$ les $\mathbb{Z}_\ell$-modules $\mathrm{H}^i(X, \mathbb{Z}_\ell) \coloneqq \lim_n \mathrm{H}^i(X, \mathbb{Z}/\ell^n)$ sont de type fini ([SGA 5 VI, 2.2]). D'après [GABBER 1983, théorème], ils sont sans torsion pour presque tout $\ell$. Si $p = 1$, cela résulte du théorème de comparaison Betti-étale d'Artin-Grothendieck ([SGA 4 XI, 4.4]) ; si $p > 1$, un ingrédient clef de sa démonstration est le *théorème du pgcd* ([DELIGNE 1980, §4.5]).

L'objet de ce paragraphe est d'expliquer une approche alternative — qui nous a été communiquée par N. Katz — au théorème susmentionné d'O. Gabber (6.2.2). Signalons cependant qu'elle repose aussi, quoique indirectement, sur le théorème du pgcd.

**6.2.1.** Soient $\mathfrak{L}$ un ultrafiltre non principal sur l'ensemble des nombres premiers différents de $p$ ([Bourbaki TG, I.§4]), $\mathfrak{m}_\mathfrak{L} = \{a = (a_\ell), \{\ell : a_\ell = 0\} \in \mathfrak{L}\}$ l'idéal maximal de l'anneau $A = \prod_{\ell \neq p} \mathbb{F}_\ell$ associé et $\mathbb{F}_\mathfrak{L} = A/\mathfrak{m}_\mathfrak{L}$ son corps résiduel, de caractéristique nulle. Supposons $k$ algébrique sur un corps fini — c'est le cas crucial — et admettons un instant (voir 6.2.3) que le foncteur $X \mapsto \mathrm{H}^\star(X, \mathbb{F}_\mathfrak{L}) \coloneqq (\prod_{\ell \neq p} \mathrm{H}^\star(X, \mathbb{F}_\ell)) \otimes_A \mathbb{F}_\mathfrak{L}$ définisse une *cohomologie de Weil* ([KLEIMAN 1968, §1.2]) ; pour notre propos, finitude, dualité de Poincaré, théorème de Lefschetz faible pour les sections hyperplanes lisses et interprétation cohomologique de la fonction $\zeta$ suffisent. Comme démontré dans [KATZ et MESSING 1974, II, corollaire 1], il résulte alors de [Deligne 1980] que l'on a pour chaque $i \geqslant 0$ l'égalité

$$\dim_{\mathbb{F}_\mathfrak{L}} \mathrm{H}^i(X, \mathbb{F}_\mathfrak{L}) = b^i(X),$$

où $b^i(X)$ est le $i$-ème nombre de Betti de $X$, nombre que l'on lit par exemple sur la fonction zêta d'un « modèle » de $X$ sur un corps fini. Le groupe de cohomologie $\mathrm{H}^i(X, \mathbb{F}_\mathfrak{L})$ étant un ultraproduit des $\mathrm{H}^i(X, \mathbb{F}_\ell)$ sa dimension est la limite, suivant l'ultrafiltre $\mathfrak{L}$, de la fonction $\ell \mapsto h^i(X, \mathbb{F}_\ell) := \dim_{\mathbb{F}_\ell} \mathrm{H}^i(X, \mathbb{F}_\ell)$, bornée d'après le théorème d'uniformité propre 3.1.1. (Voir aussi *infra* pour un argument plus direct dans ce cas particulier.) En symboles :

$$\lim_{\ell, \mathfrak{L}} h^i(X, \mathbb{F}_\ell) = b^i(X).$$

Le résultat précédent étant vrai *pour tout ultrafiltre non principal* $\mathfrak{L}$, on en déduit l'égalité $h^i(X, \mathbb{F}_\ell) = b^i(X)$ pour tous les $\ell$ sauf au plus un nombre fini d'entre eux. D'autre part, on a l'égalité

$$\dim_{\mathbb{Q}_\ell} \mathrm{H}^i(X, \mathbb{Q}_\ell) = b^i(X)$$

pour tout $\ell \neq p$, de sorte que l'on a finalement $h^i(X, \mathbb{F}_\ell) = h^i(X, \mathbb{Q}_\ell)$ pour presque tout $\ell$. La suite exacte des coefficients universels,

$$0 \to \mathrm{H}^i(X, \mathbb{Z}_\ell) \otimes \mathbb{F}_\ell \to \mathrm{H}^i(X, \mathbb{F}_\ell) \to \mathrm{H}^{i+1}(X, \mathbb{Z}_\ell)[\ell] \to 0,$$

permet alors de conclure, c'est-à-dire de retrouver le théorème d'O. Gabber sur l'absence de torsion pour presque tout $\ell$. (L'injectivité de $\mathrm{H}^i(X, \mathbb{Z}_\ell) \otimes \mathbb{F}_\ell \to \mathrm{H}^i(X, \mathbb{F}_\ell)$ suffirait.)

**6.2.2**. **Théorème** (O. Gabber ; J. Suh). *Soit $X$ un schéma propre et lisse sur un corps algébriquement clos. Alors pour presque tout nombre premier $\ell$, les groupes $\mathrm{H}^i(X, \mathbb{Z}_\ell)$ sont sans torsion.*

Le cas projectif est bien connu ([Gabber 1983, théorème]) ; comme l'ont indépendamment observé O. Gabber et J. Suh ([Suh 2012, §1]), le cas propre s'y ramène.

*Démonstration (esquisse).*
*Cas projectif.* Comme expliqué ci-dessus, cela résulte immédiatement des théorèmes de Katz-Messing et d'uniformité, ce dernier sous la forme faible présentée dans le paragraphe 6.2.4 ci-dessous.
*Réduction du cas propre au cas projectif.* D'après [de Jong 1996, 4.1], il existe une altération génériquement étale $a : X' \to X$ de degré $d$ avec $X'$ projectif et lisse sur le corps de base. Fixons un entier $n$ inversible sur $X$ et notons $\Lambda = \mathbb{Z}/n\mathbb{Z}$. Le morphisme $a$ est d'intersection complète car $X$ et $X'$ sont réguliers ; on peut donc construire un morphisme $\Lambda \to a^! \Lambda$ dans la catégorie dérivée positive des faisceaux de $\Lambda$-modules, défini par la classe de Gysin ([STG XVI, 2.5.11]). Il induit par adjonction un morphisme $a_\star a^\star \Lambda = a_! a^\star \Lambda \to \Lambda$ et le composé $\Lambda \to a_\star a^\star \Lambda \to \Lambda$ est la multiplication par le degré $d$ : comme expliqué en [STG XVI, 3.5.9, démonstration], on se ramène par restriction à un ouvert dense au cas particulier élémentaire où $a$ est fini étale.

Passant à la limite sur $n = \ell^r$, on en déduit qu'il existe pour chaque $\ell$ un morphisme trace $\mathrm{R}a_\star \mathbb{Z}_\ell \to \mathbb{Z}_\ell$ tel que les morphismes composés

$$\mathrm{H}^i(X, \mathbb{Z}_\ell) \to \mathrm{H}^i(X, \mathrm{R}a_\star \mathbb{Z}_\ell) = \mathrm{H}^i(X', \mathbb{Z}_\ell) \to \mathrm{H}^i(X, \mathbb{Z}_\ell)$$

coïncident avec la multiplication par $d$ pour chaque $i$. En particulier, si $\mathrm{H}^i(X', \mathbb{Z}_\ell)$ est sans torsion et $\ell \perp d$, il en est de même de $\mathrm{H}^i(X, \mathbb{Z}_\ell)$. Ceci nous ramène au cas projectif. □

D'après [Madore et Orgogozo 2015, 0.8], on peut en principe vérifier si un $\ell$ donné est exceptionnel. Il serait intéressant de savoir borner a priori leur nombre.

**6.2.3.** Vérifions maintenant brièvement que la « cohomologie modulo $\mathfrak{L}$ », $\mathrm{H}^\star(-, \mathbb{F}_\mathfrak{L})$, est une cohomologie de Weil. Ceci est expliqué par Ivan Tomašić dans [Tomašić 2004] ; on constate que le seul point ne résultant pas immédiatement des énoncés classiques sur la cohomologie modulo $\ell$ est la finitude de ces espaces vectoriels sur $\mathbb{F}_\mathfrak{L}$. Dans [ibid., 3.2 (i), démonstration], ceci est présenté comme un *corollaire* du théorème d'O. Gabber, ce qui rendrait l'argument précédent circulaire. Cependant, notre théorème 3.1.1 en donne une démonstration indépendante.

**6.2.4.** Signalons qu'il n'est pas difficile de majorer uniformément en $\ell$ les nombres de Betti (avec ou sans support) d'un schéma algébrique sur un corps algébriquement clos $k$ (sans utiliser le théorème sur l'absence de torsion), par des arguments proches de ceux de [Katz 2001]. (Comparer également avec [Illusie 2010, §1].)

Esquissons un argument.

(i) On souhaite montrer que si $X$ est un $k$-schéma *propre*, il existe une constante $C_X$ telle que pour tout nombre premier $\ell \neq p$ et tout $\mathbb{F}_\ell$-faisceau lisse $\mathscr{L}$ sur $X$, on ait $\dim_{\mathbb{F}_\ell} \mathrm{H}^i(X, \mathscr{L}) \leqslant C_X \cdot \mathrm{rang}(\mathscr{L})$.

On se ramène par [de Jong 1996] et descente cohomologique au cas projectif lisse et on procède ensuite par récurrence. La dualité de Poincaré et le théorème de Lefschetz nous ramènent, par l'hypothèse de récurrence, à l'étude de la cohomologie en degré médian. Ou, de façon équivalente, à la majoration uniforme des caractéristiques d'Euler-Poincaré $\chi(X, \mathscr{L})$.

(ii) On montre alors que si $X$ un $k$-schéma projectif et lisse, il existe une constante $C'_X$ telle que $|\chi(X, \mathscr{L})| \leqslant C'_X \cdot \mathrm{rang}(\mathscr{L})$ pour tout nombre premier $\ell \neq p$ et tout $\mathbb{F}_\ell$-faisceau lisse $\mathscr{L}$.

On procède à nouveau par récurrence sur $\dim(X)$. Pour une courbe, c'est un corollaire de l'égalité $\chi(X, \mathscr{L}) = \mathrm{rang}(\mathscr{L}) \cdot \chi(X)$. En général, considérons un diagramme $\mathbf{P}^1_k \leftarrow X' \to X$ où $X' \to X$ est birationnel et $X' \to \mathbf{P}^1_k$ est un pinceau de Lefschetz. Par récurrence, il suffit de démontrer le résultat pour $X'$ car si $F$ est un fermé de $X$, on a $\chi(X, \mathscr{L}) = \chi_c(U, \mathscr{L}) + \chi(F, \mathscr{L})$ et de même pour $X'$. On peut donc supposer $X = X'$. D'après la formule de Grothendieck-Ogg-Šafarevič sur $\mathbf{P}^1_k$ et l'égalité $\chi(X, \mathscr{L}) = \chi(\mathbf{P}^1_k, \mathrm{R}f_\star \mathscr{L})$, il suffit de démontrer que les $R^i f_* \mathscr{L}$ sont modérés (Picard-Lefschetz), de rangs génériques uniformément bornés par une fonction linéaire en le rang (récurrence) et qu'il en est de même du rang $h^i(X_t, \mathscr{L})$ des fibres en les points singuliers $t \in \mathbf{P}^1_k$. Ceci résulte du fait que l'écart avec le rang générique est contrôlé par les cycles évanescents (de nature locale) : les fibres des $\Psi^i_f(\mathbb{F}_\ell)$ sont de rang au plus 1 et $\mathscr{L} \simeq \mathbb{F}^r_\ell$ sur les hensélisés stricts.

(iii) Montrons maintenant que si $X$ un $k$-schéma algébrique, il existe une constante $C''_X$ telle que pour tout nombre premier $\ell \neq p$ et tout entier $i$ on ait
$$\dim_{\mathbb{F}_\ell} H^i_?(X, \mathbb{F}_\ell) \leqslant C''_X,$$
où $? \in \{\emptyset, c\}$.

Par altération et descente cohomologique, on peut supposer que $X$ est le complémentaire d'un diviseur à croisements normaux dans un schéma lisse. Le cas de la cohomologie à support compact d'un tel $X$ se ramène, par récurrence sur la dimension, au cas propre et lisse déjà traité. Le cas de la cohomologie usuelle d'un schéma lisse se ramène au cas précédent par dualité de Poincaré.

**6.2.5**. **Remarque.** Signalons que par comparaison à la caractéristique nulle le théorème d'uniformité propre — sous la forme déjà établie dans [Katz et Laumon 1985] — entraîne également des résultats d'indépendance de $\ell$ en grande caractéristique. À titre d'exemple, voici un énoncé dont on trouvera une autre démonstration dans [Brünjes et Serpé 2008, 3.4].

**Proposition.** *Soient n et d des entiers. Il existe une constante $C_{n,d}$ telle que pour tout corps algébriquement clos $k$ de caractéristique $p > C_{n,d}$ et toute hypersurface $H$ de degré $d$ de $\mathbf{P}^n_k$, on ait l'égalité des nombres de Betti $b_i(H, \mathbb{Q}_\ell) = b_i(H, \mathbb{Q}_{\ell'})$ pour tout entier $i$ et toute paire de nombres premiers $\ell$ et $\ell'$ distincts de $p$.*

Il serait intéressant de savoir si l'on peut, ne serait-ce qu'en principe, calculer une telle constante $C_{n,d}$.

# Bibliographie